\documentclass[a4paper]{article}

\usepackage[margin=3cm]{geometry}
\usepackage[dvipsnames]{xcolor}
\usepackage[english]{babel}
\usepackage[pdftex]{graphicx} 
\usepackage[utf8]{inputenc}
\usepackage{mathtools} 
\usepackage{amsfonts}
\usepackage{amsthm} 
\usepackage[T1]{fontenc}
\usepackage[pdfpagelabels]{hyperref} 
\usepackage{verbatim}
\usepackage{algorithm}
\usepackage{framed}
\usepackage{booktabs,array}

\newcommand{\R}{\mathbb{R}}
\newcommand{\lk}{\left\{}
\newcommand{\rk}{\right\}}
\newcommand{\defas}{\coloneqq}
\newcommand{\g}{\gamma}

\newcommand{\Lu}{\mathcal{L}}
\newcommand{\OuO}{\Omega \cup \Omega_I}
\newcommand{\hk}{{\varphi_k}}
\newcommand{\hj}{{\varphi_j}}
\newcommand{\sing}{\texttt{sing}}
\newcommand{\dis}{\texttt{dis}}
\newcommand{\rad}{\texttt{rad}}
\newcommand{\idx}{\texttt{idx}}
\newcommand{\X}{\texttt{X}}
\DeclareMathOperator{\bi}{\mathbf{i}}
\DeclareMathOperator{\bj}{\mathbf{j}}
\DeclareMathOperator{\spann}{span}
\DeclareMathOperator{\bL}{\mathbf{L}}
\DeclareMathOperator{\bN}{\mathbf{N}}
\DeclareMathOperator{\bn}{\mathbf{n}}
\DeclareMathOperator{\bb}{b}

\theoremstyle{plain}
\newtheorem{theo}{Theorem}[section]

\theoremstyle{definition}

\theoremstyle{remark}
\newtheorem{re}{Remark}

\newtheorem{prop}{Proposition}

\author{ Christian Vollmann\thanks{Universitaet Trier, D-54286 Trier, Germany, Email:  vollmann@uni-trier.de, volker.schulz@uni-trier.de} \and Volker Schulz\footnotemark[1] }

\title{\textbf{Exploiting multilevel Toeplitz structures in high dimensional nonlocal diffusion}}
\date{}

\begin{document}
	\maketitle
	
	\small
	\textbf{Abstract.} We present a finite element implementation for the steady-state nonlocal Dirichlet problem with homogeneous volume constraints. Here, the nonlocal diffusion operator is defined as integral operator characterized by a certain kernel function. We assume that the domain is an arbitrary $d$-dimensional hyperrectangle and the kernel is translation invariant. Under these assumptions, we carefully analyze the structure of the stiffness matrix resulting from a continuous Galerkin method with multilinear elements and exploit this structure in order to cope with the curse of dimensionality associated to nonlocal problems. For the purpose of illustration we choose a particular kernel, which is related to space-fractional diffusion and present numerical results in 1d, 2d and for the first time also in 3d.
	\normalsize

	\section{Introduction}
	The field of nonlocal operators attracts increasing attention from the mathematical society. This is due to the steadily growing pool of applications where nonlocal models are in use; including e.g., image processing \cite{imageproc1,imageproc2}, machine learning \cite{machinelearning}, peridynamics \cite{peridym1,peridym2}, fractional diffusion \cite{fractlapl} or nonlocal Dirichlet Forms \cite{nonlocaldirichletform} and jump processes \cite{jumpprocess}.
	
	In contrast to local diffusion problems, interactions can occur at distance in the nonlocal case. This relies on the definition of the nonlocal diffusion operator $-\Lu$, which acts on a function ${u\colon \R^d \to \R}$ by
	$$ -\Lu u(x) \defas 2\int_{\R^d} (u(x)-u(y))\g(x,y)dy, $$
	where $\g \colon \R^d \times \R^d \to \R$ is a nonnegative and symmetric function characterizing the precise nonlocal diffusion. In this paper we are interested in the \emph{steady-state nonlocal Dirichlet problem with volume constraints} given by
	\begin{align}
	-\Lu u(x) = f(x) ~~~&(x\in \Omega) ,\nonumber\\
	u(x) = g(x) 	~~~&(x\in \Omega_I), \label{prob}
	\end{align}
	where $\Omega \subset \R^d$ is a bounded domain. Here, the constraints are defined on a volume $\Omega_I$, the so called interaction domain, which is disjoint from $\Omega$.
	
	The recently developed vector calculus by Gunzburger et al. \cite{nonlocalveccal} builds a theoretic foundation for the description of these nonlocal diffusion phenomena. In particular, this framework allows us to consider finite-dimensional approximations using Galerkin methods similar to the analysis of (local) partial differential equations. However, in contrast to local finite element problems, here we are faced with two basic difficulties. On the one hand, the assembling procedure may require sophisticated numerical integration tools in order to cope with possible singularities of the kernel function. 
	On the other hand, discretizing nonlocal problems leads to densely populated systems. The latter tremendously affects the solving procedure, especially in higher dimensions. Thus, numerical implementations are challenging and in order to lift the concept of nonlocal diffusion from a theoretical standpoint to an applicable approach in practice, the development of efficient algorithms, which go beyond preliminary cases in 1d and 2d, is essential in this context. 
	
	In recent works, several approaches for discretizing problem (\ref{prob}) have been presented. We want to mention for instance a 1d finite element code by D'Elia and Gunzburger \cite{fractlapl}, where a truncated version of the fractional kernel $\g(x,y) =\frac{c_{d,s}}{2||y-x||^{d+2s}}$ has been used, but which can easily be extended to general (singular) kernels. Also for the fractional kernel, Acosta, Bersetche and Borthagaray \cite{acosta} developed a finite element implementation for the 2d case. In general, a lot of work has been done for the discretization of fractional diffusion problems and fractional derivatives in various definitions (mainly via finite difference schemes) not only for 1d and 2d \cite{wangToeplitz,wang2d}, but also for the 3d case \cite{wang3d}. However, to the best of out knowledge, for general kernel functions 3d finite element implementations for problem (\ref{prob}) are not yet available.
	
	In this paper we study a finite element approximation for problem (\ref{prob}) on an arbitrary $d$-dimensional hyperrectangle (parallel to the axis) for translation and reflection invariant kernel functions. More precisely, we analyze from a computational point of view a continuous Galerkin discretization with multilinear $Q_1$ elements for the following setting:
	\begin{itemize}
		\item[(A1)] We set $\Omega \defas \prod_{i=0}^{d-1}[a_i,b_i]$, where $[a_i, b_i]$ are compact intervals on $\R$.
		\item[(A2)] We assume that the kernel $\g$ is \emph{translation and reflection invariant}, such that 
		$$\g(b+R_ix, b+R_iy) = \g(x,y)$$
		for all $b\in\R^d$ and all $0\leq i \leq d$, where $R_i(x) \defas (x_0,\ldots,-x_i,\ldots,x_{d-1})$ and $R_d \defas Id$. 
	\end{itemize}
	As a consequence, these structural assumptions on the underlying problem are reflected in the stiffness matrix; we obtain a symmetric $d$-level Toeplitz matrix, which has two crucial advantages. On the one hand, we only need to assemble (and store) the first row (or column) of the stiffness matrix. On the other hand, we can benefit from an efficient implementation of the matrix-vector product for solving the linear system. 
	This result is presented in Theorem \ref{multilevel_toeplitz} and is crucial for this work, since it finally enables us to solve the discretized system in an affordable way. 
	For illustrative purposes we choose the fractional kernel and exploit a third assumption on the interaction horizon for simplifying the implementation:
	\begin{itemize}
		\item[(A3)] We assume that interactions only occur at a certain distance $R$, which we assume to be larger or equal the diameter of the domain $\Omega$, such that $\Omega \subset B_R(x)$ for all $x \in \Omega$.
	\end{itemize}	
	
	The paper is organized as follows. In Section \ref{sec:theory} we cite the basic results about existence and uniqueness of weak solutions and finite-dimensional approximations. In Section \ref{sec:femsetting} we give details about the precise finite element setting and proof our main result, that the stiffness matrix is multilevel Toeplitz. In Sections \ref{sec:assembling} and \ref{sec:solving} we explain in detail the implementation of the assembling and solving procedure, respectively. In Section \ref{sec:numres} we round off these considerations by presenting numerical results with application to space-fractional diffusion.

	\section{Nonlocal diffusion problems}\label{sec:theory}
	We review the relevant aspects of nonlocal diffusion problems as they are introduced in \cite{wellposedness} which constitute the theoretic fundamentals of this work. 
	
	Let $\Omega \subset \R^d$ be a bounded domain with piecewise smooth boundary. Further let $\g \colon \R^d \times \R^d \to \R$ be a nonnegative and symmetric function (i.e., $\g(x,y) = \g(y,x) \geq 0$), which we refer to as \emph{kernel}. Then we define the action of the \textit{nonlocal diffusion operator} $-\Lu \defas -\Lu_\g $ on a function $u \colon \R^d \to \R$ by 
	$$ -\Lu u(x) \defas 2\int_{\R^d} (u(x)-u(y))\g(x,y)dy $$
	for $x \in \Omega$. In addition to that, we assume that there exists a constant $\g_0 > 0$ and a finite \emph{interaction horizon} or \emph{radius} $R > 0$ such that for all $x \in \Omega$ we have
	\begin{align}
	\g(x,y) \geq 0 ~~~&\forall~ y \in B_R(x), \nonumber\\
	\g(x,y) \geq \g_0 > 0 ~~~&\forall ~y \in B_{R/2}(x),\nonumber\\
	\g(x,y) = 0   ~~~&\forall ~y \in (B_{R}(x))^c, \label{gen_assumption}
	\end{align}
	where $B_R(x) \defas \lk y \in \R^d \colon ||x-y||_2< R\rk$. Thus we can consider the kernel as a composition
	\begin{align}
	\g(x,y) =g(x,y)\mathcal{X}_{B_R(x)}(y),~~x,y \in \R^d , \label{kernel_indicatorfunction}
	\end{align}
	for some appropriate nonnegative and symmetric function $g$ and the indicator function $\mathcal{X}_{B_R(x)}(y)$. Unless otherwise stated, we from now on consider the indicator function as part of the kernel, such that we can notationally omit the intersection between the ball $B_R(x)$ and the domain of integration for integrals involving the kernel. Further, we define the \emph{interaction domain} by $$\Omega_I \defas \lk y \in \Omega^c\colon \exists x\in \Omega\colon\g(x,y) \neq 0 \rk$$ and finally introduce the \emph{steady-state nonlocal Dirichlet problem with volume constraints} as
	\begin{align*}
	-\Lu u(x) = f(x) ~~~&(x\in \Omega) ,\nonumber\\
	u(x) = g(x) 	~~~&(x\in \Omega_I), 
	\end{align*}
	where $f\colon \Omega \to \R$ is called the \emph{source} and $g: \Omega_I \to \R$ specifies the Dirichlet volume constraints. In the remainder of this paper we assume $g\equiv 0$.
	\subsection{Weak formulation}
	For the purpose of constructing a finite element framework for nonlocal diffusion problems, we introduce the concept of weak solutions as it is presented in \cite{wellposedness}. 
	
	We define the bilinear form
	\begin{align*}
	a(u,v) \defas& \int_{\Omega} v(-\Lu u )  dx,
	\end{align*}
	and the associated linear functional
	\begin{align*}
	\ell(v) \defas \int_{\Omega} fvdx.
	\end{align*}
	By establishing a nonlocal vector calculus it is shown in \cite{wellposedness}, that, if $u$ and $v$ are zero on the interaction domain, the following equality holds:
	\begin{align*}
	a(u,v)= &\int_{\OuO} \int_{\OuO} (u(y)-u(x))(v(y)-v(x))\g(x,y)  dydx.
	\end{align*}
	This implies that $a$ is symmetric and nonnegative, or equivalently, the linear nonlocal diffusion operator $-\Lu$ is self-adjoint with respect to the $L^2$-product and nonnegative. Furthermore, we define the \emph{nonlocal energy space} 
	$$V(\OuO) \defas \lk u \in L^2(\OuO) \colon |||u||| <\infty \rk,$$
	where 
	$ |||u|||\defas \sqrt{\tfrac{1}{2}a(u,u)}$
	and the \emph{nonlocal constrained energy space} 
	$$V_c(\OuO) \defas \lk u \in V(\OuO) \colon  u_{|\Omega_I} \equiv 0 ~\text{a.e.}\rk.$$	
	We note that $||| \cdot |||$ constitutes a semi-norm on $V(\OuO) $ and due to the volume constraints a norm on $V_c(\OuO) $. With these preparations at hand, a weak formulation of (\ref{prob}) can be formulated as
	\begin{align}
	\text{\textit{Find $u \in V_c(\OuO) $ such that }} a(u,\cdot) \equiv \ell(\cdot) ~ \text{\textit{on}}~ V_c(\OuO) \label{weakform}.
	\end{align}
	
	In order to make statements about the existence and uniqueness of weak solutions we have to further specify the kernel. In \cite{wellposedness} the authors consider, among others, a certain class of kernel functions, on which we will focus in the remainder of this section and in our numerical experiments. More precisely, we require that there exists a fraction $s\in (0,1)$ and constants $\g_1,\g_2 > 0$, such that for all $x \in \OuO$ it holds that
	\begin{align}
	\g_1\leq \g(x,y)||y-x||_2^{d+2s} \leq\g_2 ~~~&\forall ~y \in B_R(x) \label{case_1}.
	\end{align}
	Then it is shown in \cite{wellposedness} that the nonlocal constrained energy space $(V_c(\OuO),|||\cdot|||)$ is equivalent to the constrained fractional-order Sobolev space 
	\begin{align}
	H^s_c(\OuO) \defas \lk u \in H^s(\OuO)\colon u_{|\Omega_I} \equiv 0 ~\text{a.e.} \rk,
	\end{align}
	where
	\begin{align}
	H^s(\OuO) \defas \lk u \in L^2(\OuO) \colon ||u||_{H^s(\OuO)} \defas ||u||_{L^2(\OuO)} +  |u|_{H^s(\OuO)} < \infty   \rk
	\end{align}
	and 
	$$|u|_{H^s(\OuO)}^2 \defas  \int_{\OuO} \int_{\OuO}  \frac{(u(x)-u(y))^2}{||x-y||_2^{d+2s}}dy dx.   $$
	Hence, there exist two positive constants $C_1$ and $C_2$ such that
	\begin{align}
	C_1 ||u||_{H^s(\OuO)} \leq |||u||| \leq C_2 ||u||_{H^s(\OuO)} ~~~\forall~ u \in V_c(\OuO). \label{space_eq}
	\end{align}
	This equivalence implies that $(V_c(\OuO),|||\cdot|||)$ is a Banach space for kernel functions satisfying (\ref{gen_assumption}) and (\ref{case_1}). Applying Lax-Milgram Theorem finally brings in the well posedness of problem (\ref{weakform}); see \cite{wellposedness}.
	
	\subsection{Finite-dimensional approximation}
	With the concept of weak solutions we can proceed as in the local case to develop finite element approximations of (\ref{weakform}). 
	
	Therefore, let $\lk V_c^N\rk_N$ be a sequence of finite-dimensional subspaces of $V_c(\OuO)$, where $N = \dim (V_c^N)$, and let $u_N$ denote the solution of
	\begin{align}
	\emph{\text{Find $u_N \in V_c^N$ such that } $a(u_N, v_N) = \ell(v_N)$ \text{ for all $v_N$ in $V_c^N$} }. \label{finiteelementproblem}
	\end{align}
	Then from \cite{Olena} we recall the following regularity and convergence results.
	\begin{prop}[{\cite[Theorem 3.5, Proposition 3.6]{Olena}}]
		Let the domain $\Omega \subset \R^d$ have $C^\infty$ boundary $\partial \Omega$ and let $f \in H^r(\Omega)$ for $r\geq 0$. Further let the kernel be of the form
		\begin{align*}
		\g(x,y) = \frac{c}{||x-y||_2^{d+2s}} \mathcal{X}_{B_R(x)}(y),
		\end{align*} 
		for a constant $c > 0$, such that (\ref{gen_assumption}) and (\ref{case_1}) are satisfied. Then for the solution $u \in V_c(\OuO)$ of (\ref{weakform}) the following regularity estimate holds
		\begin{align*}
		|u|_{H^{s+\alpha}(\OuO)} \leq C ||f||_{H^r(\Omega)},~~~C>0,
		\end{align*}
		where $\alpha = \min \lk s+r,1/2-\varepsilon \rk $ for some arbitrarily small $\varepsilon >0$. Furthermore, by invoking this regularity estimate we obtain the following convergence result for piecewise linear finite element approximations:
		\begin{align}
		||u-u^h||_{H^s(\OuO)} \leq C' h^\alpha ||f||_{H^r(\Omega)},~~~C'>0. \label{convergence_rate}
		\end{align}
	\end{prop}
	
	To the best of our knowledge, for the truncated fractional kernel and less smooth domains corresponding results are not available. However, for Lipschitz domains and the untruncated fractional kernel
	\begin{align*}
	\g_\infty(x,y) = \frac{c}{||x-y||_2^{d+2s}}
	\end{align*}
	similar regularity estimates have been obtained under the Hölder regularity assumption on the right-hand side $f$ \cite{borthagaray_regularity}.\\
	
	The derivation of the discretized problem then relies on the construction of the stiffness matrix. Therefore let $\lk \varphi_0, \ldots, \varphi_{N-1} \rk$ be a basis of $V_c^N$, such that the finite element solution $u^N \in V_c^N$ can be expressed as a linear combination 
	$u^N = \sum_{k=0}^{N-1} u_k^N \hk $. If the basis functions are chosen in a way such that $\hk(x_k) = 1$ on appropriate grid points $x_k$, then the coefficients satisfy $u_k^N = u^N(x_k)$. We test for all basis functions, such that the finite element problem (\ref{finiteelementproblem}) reads as
	\begin{align*}
	\textit{Find $u^N\in \R^{N}$ such that }\sum_{k=0}^{N-1} u_k^N a(\hk, \hj) = \ell(\hj) =: b_j  \textit{ for } 0 \leq j < N.
	\end{align*}
	The stiffness matrix $A^N=(a_{kj})_{kj} \in \R^{N \times N}$ is given by
	\begin{align*}
	a_{kj} \defas  \int_{\OuO} \int_{\OuO}  (\hk(y)-\hk(x))(\hj(y)-\hj(x))\g(x,y)dydx
	\end{align*}
	and we finally want to solve the discretized \emph{Galerkin system}
	\begin{align}
	A^Nu^N = b^N, \label{discretizedproblem}
	\end{align}
	where $u^N, b^N \in \R^{N}$. The properties of the bilinear form $a$ imply that $A^N$ is symmetric and positive definite, such that there exists a unique solution $u^N$ of the finite-dimensional problem (\ref{discretizedproblem}).
	
	\section{Finite element setting} \label{sec:femsetting}
	In this section we study a continuous Galerkin discretization of the homogeneous nonlocal Dirichlet problem, given by
	\begin{align*}
	2\int_{\OuO} (u(x)-u(y))\g(x,y)dy = f(x) ~~~&(x\in \Omega) ,\nonumber\\
	u(x) = 0 	~~~~~~~&(x\in \Omega_I),
	\end{align*}
	under the following assumptions:
	\begin{itemize}
		\item[(A1)] We set $\Omega \defas \prod_{i=0}^{d-1}[a_i,b_i]$, where $[a_i, b_i]$ are compact intervals on $\R$.
		\item[(A2)] We assume that the kernel $\g$ is \emph{translation and reflection invariant}, such that 
		$$\g(b+R_ix, b+R_iy) = \g(x,y)$$
		for all $b\in\R^d$ and all $0\leq i \leq d$, where $R_i(x) \defas (x_0,\ldots,-x_i,\ldots,x_{d-1})$ and $R_d \defas Id$.
	\end{itemize}
	Assumption (A1) allows for a simple triangulation of $\Omega$, which we use to define a finite-dimensional energy space $V^{N}_c$. Together with (A2) we can show that this discretization yields the multilevel Toeplitz structure of the stiffness matrix, where the order of the matrix is determined by the number of grid points in each respective space dimension. 
	\subsection{Definition of the finite-dimensional energy space}
\begin{figure}
		\centering
		\def\svgwidth{\textwidth}
		\begingroup%
		\makeatletter%
		\providecommand\color[2][]{%
			\errmessage{(Inkscape) Color is used for the text in Inkscape, but the package 'color.sty' is not loaded}%
			\renewcommand\color[2][]{}%
		}%
		\providecommand\transparent[1]{%
			\errmessage{(Inkscape) Transparency is used (non-zero) for the text in Inkscape, but the package 'transparent.sty' is not loaded}%
			\renewcommand\transparent[1]{}%
		}%
		\providecommand\rotatebox[2]{#2}%
		\ifx\svgwidth\undefined%
		\setlength{\unitlength}{183.94099919bp}%
		\ifx\svgscale\undefined%
		\relax%
		\else%
		\setlength{\unitlength}{\unitlength * \real{\svgscale}}%
		\fi%
		\else%
		\setlength{\unitlength}{\svgwidth}%
		\fi%
		\global\let\svgwidth\undefined%
		\global\let\svgscale\undefined%
		\makeatother%
		\begin{picture}(1,0.32205233)%
		\put(0,0){\includegraphics[width=\unitlength,page=1]{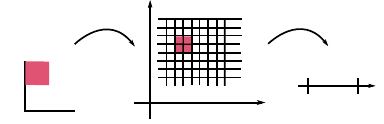}}%
		\put(0.79135697,0.03804194){\color[rgb]{0,0,0}\makebox(0,0)[lb]{\smash{$0$}}}%
		\put(0.92491334,0.03791293){\color[rgb]{0,0,0}\makebox(0,0)[lb]{\smash{$1$}}}%
		\put(0.60837399,0.01299861){\color[rgb]{0,0,0}\makebox(0,0)[lb]{\smash{$b_0$}}}%
		\put(0.33602039,0.26366791){\color[rgb]{0,0,0}\makebox(0,0)[lb]{\smash{$b_1$}}}%
		\put(0.20497551,0.06513067){\color[rgb]{0,0,0}\makebox(0,0)[lb]{\smash{$1$}}}%
		\put(0.24814632,0.26261517){\color[rgb]{0,0,0}\makebox(0,0)[lb]{\smash{$T_k$}}}%
		\put(0.76612339,0.27121812){\color[rgb]{0,0,0}\makebox(0,0)[lb]{\smash{$\varphi_k$}}}%
		\put(0.40395794,0.0121393){\color[rgb]{0,0,0}\makebox(0,0)[lb]{\smash{$a_0$}}}%
		\put(0.55630574,0.297979){\color[rgb]{0,0,0}\makebox(0,0)[lb]{\smash{$x_k$}}}%
		\put(0,0){\includegraphics[width=\unitlength,page=2]{basisfunction.pdf}}%
		\put(-0.00438558,0.06619341){\color[rgb]{0,0,0}\makebox(0,0)[lb]{\smash{$-1$}}}%
		\put(0,0){\includegraphics[width=\unitlength,page=3]{basisfunction.pdf}}%
		\put(0.32841117,0.08977666){\color[rgb]{0,0,0}\makebox(0,0)[lb]{\smash{$a_1$}}}%
		\put(0,0){\includegraphics[width=\unitlength,page=4]{basisfunction.pdf}}%
		\end{picture}%
		\endgroup%
		\caption{\small Transformation of the reference basis function, where $d=2$.}
		\label{fig:basisfunction}
	\end{figure}

	We decompose the domain $\Omega = \prod_{i=0}^{d-1}[a_i,b_i]$ into $d$-dimensional hypercubes with sides of length $ h>0 $ in each respective dimension. Note that we can omit a discretization of $\Omega_I$ since we assume homogeneous Dirichlet volume constraints. 
	Let $\bN = (N_i)_{0 \leq i < d} \defas (\frac{b_i-a_i}{ h })_{0 \leq i < d}$ and ${\bL \defas (N_i-1)_{0 \leq i < d}}$, then for the interior of $\Omega$ this procedure results in $\bL^d \defas \prod_{i=0}^{d-1} L_i$ degrees of freedom. Due to the simple structure of the domain we can choose a canonical numeration for the resulting grid $\prod_{i=0}^{d-1} \left( a_i +  h \lk 0, \ldots, L_i -1\rk \right) $ of inner points. More precisely, we will employ the map
	\begin{align*}
	E^{\bn}(z) \defas \sum_{i=0}^{d-1} z_i p_i(\bn),
	\end{align*}
	where $p_i(\bn) \defas \prod_{j>i}n_j$, for establishing an order on a structured grid $\prod_{i=0}^{d-1} \lk 0, \ldots n_i-1 \rk$, where $\bn=(n_0,\ldots,n_{d-1}) \in \mathbb{N}^d$.
	Its inverse is given by
	\begin{align*}
	\left(E^{-\bn}(k)\right)_{0 \leq i < d} =  \left(\lfloor \tfrac{k}{p_i(n)}\rfloor - \lfloor \tfrac{k}{p_{i-1(n)}}\rfloor n_i \right)_{0 \leq i < d}.
	\end{align*}
	Let $e \defas (1,\ldots, 1)\in \R^{d}$ and $a \defas (a_0,\ldots,a_{d-1})$, then we define the ordered array of inner grid points
	$
	(x_k)_{0 \leq k < \bL^d} \in \R^{\bL^d \times d} 
	$
	by $$ x_k   \defas  a  + h ( E^{-\bL}(k)+e)$$ for $0 \leq k < \bL^d$. We further define elements $S_k \defas b_k +  h \square$, where $ \left(b_k^i\right)_{0 \leq i < d}   \defas  a +  h  E^{-\bN}(k) $ and $\square \defas [0,1]^d$, such that
	$\Omega = \bigcup_{k=0}^{\bN^d - 1} S_k.$ Next we aim to define appropriate element basis functions on the reference element $\square$. Therefore we denote by 
	$$(v_k)_{0 \leq k < 2^d} \in \R^{2^d \times d}$$
	the vertices of the unit cube $\square$ ordered according to 
	$
	v_k  \defas E^{-(2,\ldots, 2)}(k).
	$
	Then for each vertex $v_k$, $0 \leq k < 2^d$, we define an \emph{element basis function} $\psi_k \colon \square \to [0,1]$ by
	\begin{align*}
	\psi_k(x) &= \left(\prod_{i=0, v_k^i = 0}^{d-1} (1-x_i) \right)\left(\prod_{i=0, v_k^i = 1}^{d-1} x_i \right)	.
	\end{align*}
	For dimensions $d \in \lk 1,2,3 \rk$ respectively, these are the usual linear, bilinear and trilinear element basis functions (see e.g. \cite[Chapter 1]{fem_book}). They are defined in a way such that $0 \leq \psi_k \leq 1$ and $\psi_k(v_k) = 1$. Moreover, we define the \emph{reference basis function} $\varphi \colon \R^d \to [0,1]$ by
	\begin{align}
	\varphi(x) &\defas \begin{cases}
	\psi_i(v_i + x)&: x \in (\square - v_i)\\
	0&:else.
	\end{cases} \label{reference_basisfunction}
	\end{align}
	We note that $J\defas [-1,1]^d = \dot\bigcup_{i=0}^{2^d} (\square -v_i ) $ (disjoint union), such that $\varphi$ is well defined and $supp(\varphi) =  J$. Now let the physical support be defined as
	$$I_k \defas \bigcup \lk S_i \colon  x_k \in S_i, 0 \leq i < \bN^d \rk ,$$
	which is a patch of the elements touching the node $x_k$. We associate to each element $S_k$ the transformation $T_k \colon  J \to I_k$, $T_k(v) \defas x_k +  h  v $. We note that $\det dT_k(x) \equiv h^d  $.
	Then for each node $x_k$ we define a basis function $\varphi_k \colon \OuO \to [0,1]$ by
	\begin{align*}
	\varphi_k(x) &\defas \begin{cases}
	\varphi(T_k^{-1}(x))&: x \in I_k\\
	0&:else
	\end{cases}\\
	&=\begin{cases}
	\psi_i(v_i + T_k^{-1}(x))&: T_k^{-1}(x)\in (\square - v_i)\\
	0&:else,
	\end{cases}
	\end{align*}
	which satisfies $0\leq \varphi_k \leq 1$ and $\varphi_k(x_k) = 1$. Figure \ref{fig:basisfunction} illustrates the latter considerations for $d=2$.
	Finally, we can define a constrained finite element space by
	\begin{align}
	V_c^{ h } \defas \spann \lk \varphi_k \colon 0 \leq k < \bL^d \rk \label{fem_space},
	\end{align}
	such that each linear combination consisting of a set of these basis functions fulfills the homogeneous Dirichlet volume constraints. Notice that we parametrize these spaces by the grid size $ h $ indicating the dimension $\bL^d$, which is by definition a function of $ h $. Finally, we close this subsection with the following observations, which we exploit in the remainder.
	\begin{re}
		Let $x \in \R^d$, then:
		\begin{itemize}
			\item[i)]  $\varphi(x) = \varphi(|x|)$, where $|x| \defas ( |x_i|)_i$.
			\item[ii)] Let $R_i \colon \R^d \to \R^d$, for $0 \leq i < d$, denote the reflection 
			$
			R_i(x) = (x_0, \ldots, -x_i, \ldots, x_{d-1}),
			$ then		
			\begin{align}
			\varphi(x)  = \varphi (|x| ) \Leftrightarrow \varphi(x)  = \varphi (R_i(x)) ~~\forall~ 0 \leq i < d. \label{reflection_invariance}
			\end{align}
			\item[iii)]  $\varphi(x)	= \varphi((x_{\sigma(i)})_i)$ for all permutations $\sigma \colon \lk 0,\ldots, d-1 \rk \to \lk 0,\ldots, d-1 \rk$.
		\end{itemize}

	\end{re}
	\begin{proof} We first show i).
		Since $\varphi(x) = 0 = \varphi(|x|)$ for $x$ in $int({J})^c$, let $x \in int({J})$. Thus, there exists an index $0 \leq k < d $ such that 
		$x \in int(\square) - v_k,$ which implies that 
		$x_i < 0$ if and only if $v_k^i = 1. $ Hence, we can conclude that
		\begin{align*}
		\varphi(x) &= \psi_k(x + v_k) \\
		&= \left(\prod_{i=0, v_k^i = 0}^{d-1} (1-x_i) \right)\left(\prod_{j=0, v_k^i = 1}^{d-1} (1+x_i) \right) \\
		&= \left(\prod_{i=0, v_k^i = 0}^{d-1} (1-|x_i|) \right)\left(\prod_{j=0, v_k^i = 1}^{d-1} (1-|x_i|) \right) \\
		&= \prod_{i=0}^{d-1} (1-|x_i|) \\
		&=\psi_0 (|x| + v_0) \\
		&= \varphi(|x|).
		\end{align*}
		Then, on the one hand, we have that $|x| = |R_i(x)|$ and therefore
		$\varphi(x)  = \varphi (R_i(x))$ for all $0 \leq i < d$.
		On the other hand, we note that the operation $|\cdot|$ is a composition of reflections $R_i$, more precisely
		$$|x| = \left(\prod_{x_i < 0}R_i\right)(x) .$$ Thus, we obtain the equivalence stated in ii). Statement iii) follows from the representation
		$
		\varphi(x)	= \varphi(|x|) 
		=\prod_{i=0}^{d-1} (1-|x_i|)
		$
		due to the commutativity of the product.
	\end{proof}

	\subsection{Multilevel Toeplitz structure of the stiffness matrix}
	Now we aim to show that the stiffness matrix $A$ owns the structure of a $d$-level Toeplitz matrix. This is decisive for this work, since it finally enables us to solve the discretized system (\ref{discretizedproblem}) in an affordable way. 
	
	From now on the assumption (A2) on the kernel function becomes crucial. At this point we note that the indicator function
	$(x,y) \mapsto \mathcal{X}_{B_R(x)}(y) $
	is translation and reflection invariant. This even holds if the ball $B_R(x)$ is defined with respect to another than the $||\cdot||_2$-norm. Hence, as mentioned in \eqref{kernel_indicatorfunction}, we can regard the kernel as a composition
	$$ \g(x,y) = g(x,y)\mathcal{X}_{B_R(x)}(y),~~x,y \in \R^d ,$$
	for some translation and reflection invariant function $g$. In order to analyze the multilevel structure of ${A \in \R^{\bL^d \times \bL^d}}$ it is convenient to introduce an appropriate multi-index notation. 
	To this end, we choose $E^{\bL}$ from above as index bijection
	and we identify
	$a_{\bi\bj} = a_{E^{\bL}(\bi),E^{\bL}(\bj)}. $
	We call the matrix $A$ \emph{$d$-level Toeplitz} if 
	$$
	a_{\bi\bj} = a(\bi -\bj).
	$$
	If even $a_{\bi\bj} = a(|\bi -\bj|)$, where the absolute value is understood componentwise, then each level is symmetric and we can reconstruct the whole matrix from the first row (or column). For a more general and detailed consideration of multilevel Toeplitz matrices see for example \cite{multileveltoeplitz}. However, with this notation at hand we can now formulate
	\begin{theo} \label{multilevel_toeplitz}
		Let the kernel $\g$ fulfill assumptions (A1) and (A2) and let the finite element space $V_c^{ h }$ be defined as in (\ref{fem_space}) for a grid size $ h >0$. Then the stiffness matrix $A$ associated with problem (\ref{finiteelementproblem}) is $d$-level Toeplitz, where each level is symmetric. 
	\end{theo}
	\begin{proof}
		The key point in the proof is the relation
		\begin{align*}
		a_{kj} = a\left(\left|{ h }^{-1}(x_k-x_j)\right|\right),
		\end{align*}
		which we show in two steps. First we show that $a_{kj} = a\left( { h }^{-1}(x_k-x_j)\right) $ and then we proof $a(z) = a(|z|).$ Therefore let us recall that the entry $a_{kj}$ of the stiffness matrix $A$ is given by 
		$$a_{kj} = \int_{\OuO} \int_{\OuO }  (\hk(y)-\hk(x))(\hj(y)-\hj(x))\g(x,y)dydx.$$ 
		Having a closer look at the support of the integrand, we find that 
		\begin{align*}
		(\hk(y)-\hk(x))(\hj(y)-\hj(x)) = 0  \Leftrightarrow (x,y) \in (I_k^c \times I_k^c) \cup (I_j^c \times I_j^c) \cup \lk (x,x)\colon x \in \OuO \rk.
		\end{align*} 
		Since $\lk (x,x)\colon x \in \OuO\rk$ has null $\lambda_{2d}$-Lebesgue measure we can neglect it in the integral and obtain
		\begin{align*}
		a_{kj} = \int_{\left( I_k^c \times I_k^c\right)^c \cap \left(I_j^c \times I_j^c \right)^c} (\hk(y)-\hk(x))(\hj(y)-\hj(x))\g(x,y)dydx.
		\end{align*}
		\\
		Aiming to show $a_{kj} = a\left({ h }^{-1}(x_k-x_j)\right)$ we need to carry out some basic transformations of this integral. Since by definition $\varphi_j = \varphi \circ T_j^{-1}$ and also $\det T_j(x) \equiv h^d$ we find
		\begin{align*}
		a_{kj} &= \int_{\left( I_k^c \times I_k^c\right)^c \cap \left(I_j^c \times I_j^c \right)^c} 
		(\hk(y)-\hk(x))(\hj(y)-\hj(x))\g(x,y)dydx\\
		&=  h ^{2d}\int_{T_j^{-1}(\R^d) \times T_j^{-1}(\R^d) }
		\left(1-\mathcal{X}_{I_j^c \times I_j^c }(T_j(v),T_j(w))\right)
		\left(1-\mathcal{X}_{I_k^c \times I_k^c }(T_j(v),T_j(w))\right)\\
		&~~~~~~~~~~~~ \left((\varphi \circ  T_k^{-1})(T_j(w))-(\varphi \circ T_k^{-1})(T_j(v))\right)\left(\varphi(w)-\varphi(v)\right)\g\left(T_j(v),T_j(w)\right)dwdv.
		\end{align*}
		Now we make a collection of observations. Due to assumption (A2) we have
		$$\g(T_j(v),T_j(w)) = \g(x_j +  h  v,x_j +  h  w) = \g( h  v,  h  w). $$ Furthermore, by definition of the transformations $T_j, T_k$ we find that $T_j^{-1}(\R^d)=\R^d$ as well as
		$(T_k^{-1}\circ T_j)(v) =   h ^{-1}(x_j +  h  v - x_k) =  h ^{-1}(x_j - x_k)+ v $. Since these transformations are bijective we also have that
		$\mathcal{X}_{M^c\times M^c }(T_j(x),T_j(y)) = \mathcal{X}_{(T_j^{-1}(M))^c\times (T_j^{-1}(M))^c }(x,y)$ for a set $M\subset \R^d$. Hence, defining $x_{jk}  \defas  h ^{-1}(x_j - x_k) = - x_{kj} $ and recognizing $T_j^{-1}(I_k) = x_{kj} + J$ we finally obtain
		\begin{align*}
		a_{kj}
		&= { h }^{2d}\int_{\left( J^c \times J^c\right)^c \cap \left((x_{kj}+J )^c \times (x_{kj}+J )^c\right)^c}
		\\
		&~~~~~~~~~~~~ \left(\varphi(w-x_{kj})-\varphi(v-x_{kj})\right)\left(\varphi(w)-\varphi(v)\right)\g( h  w, h  v)dwdv\\
		&=a(x_{kj}).
		\end{align*} 
		Next, we proof that this functional relation fulfills $a(z) = a(|z|)$.
		Let us for this purpose define 
		$F(x,y;z) \defas  (\varphi(y-z)-\varphi(x-z))(\varphi(y)-\varphi(x))\g( h  y, h  x)$
		such that
		$$a(z) = { h }^{2d} \int_{\left( J^c \times J^c\right)^c \cap \left((z+J )^c \times (z+J )^c \right)^c}
		F(x,y;z)dydx. $$
		Let $z \in \lk x_{kj}: 0 \leq k,j<\bL^d \rk$. Then there exists a matrix $R= R(z)  \in \R^{d \times d}$, which is a composition of reflections $R_i$ from (\ref{reflection_invariance}), such that $Rz = |z|$. Then from (\ref{reflection_invariance}) and the assumption (A2) on the kernel, we obtain for $x,y \in \R^d$ that
		\begin{align*}
		F(Rx,Ry;|z|) & = \left(\varphi(Ry-Rz)-\varphi(Rx-Rz)\right)\left(\varphi(Ry)-\varphi(Rx)\right)\g( h  Ry, h  Rx) \\
		&= \left(\varphi(y-z)-\varphi(x-z)\right)\left(\varphi(y)-\varphi(x)\right)\g( h  y, h  x) \\
		&=F(x,y;z).
		\end{align*}
		Since
		$R(J ) = J  $
		and therefore
		\begin{align*}
		&R\left(\left( J^c \times J^c\right)^c \cap \left((z+J )^c \times (z+J )^c\right)^c\right) \\
		=&\left( J^c \times J^c\right)^c \cap \left((|z|+J )^c \times (|z|+J )^c\right)^c,
		\end{align*}
		we eventually obtain
		\begin{align*}
		a(|z|)&={ h }^{2d} \int_{\left( J^c \times J^c\right)^c \cap \left((|z|+J )^c \times (|z|+J )^c\right)^c}
		F(x,y;|z|)dydx \\
		&={ h }^{2d} \int_{\left( J^c \times J^c\right)^c \cap \left((z+J )^c \times (z+J )^c\right)^c}
		F(Rx,Ry;|z|)dydx \\
		&={ h }^{2d} \int_{\left( J^c \times J^c\right)^c \cap \left((z+J )^c \times (z+J )^c\right)^c}
		F(x,y;z)dydx \\
		&=a(z).
		\end{align*}
		Finally, we can show that $A$ carries the structure of a $d$-level Toeplitz matrix. By having a closer look at the definitions of $E^{\bL}$ and the grid points $x_k$ we can conclude that
		\begin{align*}
		a_{\bi\bj} &= a_{E^{\bL}(\bi)E^L(\bj)}
		= a\left(\left|{ h }^{-1}(x_{E^{\bL}(\bi)}-x_{E^{\bL}(\bj)})\right|\right)
		= a(|\bi - \bj|).
		\end{align*}
		Thus, the entry $a_{\bi\bj}$ only depends on the difference $\bi - \bj$.
	\end{proof}

	With other words, the translation invariance of the kernel brings in the relation $a_{kj} = a\left({ h }^{-1}(x_k-x_j)\right)$. The advantage of this observation relies on the usage of a regular grid leading to redundancy in the set $\lk x_k-x_j \colon 0 \leq k,j < \bL^d \rk$. From the reflection invariance we can finally deduce $a_{kj} = a\left(\left|{ h }^{-1}(x_k-x_j)\right|\right)$ leading to symmetry in each level.
	
	As a consequence, in order to implement the  matrix-vector product, it is sufficient to assemble solely the first row or column
	$$M \defas (a_{\ell0})_\ell = \left(a\left({ h }^{-1}(x_\ell-x_0)\right)\right)_\ell = \left(a\left(E^{\bL}(\ell)\right)\right)_\ell$$ 
	of the stiffness matrix $A$, since for $\ell(k,j) \defas E^L\left( h ^{-1}(|x_k-x_j|)\right)$ we get 
	$$ a_{kj} = a\left( h ^{-1}\left(|x_k-x_j|\right)\right) = a\left(E^{\bL}\left(\ell(k,j)\right)\right) = M_{\ell(k,j)}. $$
	Note that $\ell(k,j) \defas E^{\bL}( h ^{-1}(|x_k-x_j|))$ is well defined, since $ h ^{-1}\left(|x_k-x_j|\right) $ lies in the domain of definition of $E^{\bL}$. 
	
	\begin{re} \label{L1_radial}
		Exploiting that $\varphi$ is invariant under permutations, the same proof (by composing the reflection $R$ with a permutation matrix) shows that
		$a(z)	= a((z_{\sigma(i)})_i)$ for all permutations $\sigma \colon \lk 0,\ldots, d-1 \rk \to \lk 0,\ldots, d-1 \rk$. We will use this observation to accelerate the assembling process. Also note in this regard, that a kernel of radial type, i.e., $\g(x,y) = g(||x-y||_2)\mathcal{X}_{B_R(x)}(y)$, is also invariant under such permutations, independent of the norm used to define the ball $B_R(x)$.
	\end{re}

	\section{Assembling procedure}\label{sec:assembling}
	In this section we aim to analyze the entries $a_{kj} = a(x_{kj})$ of the stiffness matrix $A$ more closely and derive a representation which can be efficiently implemented. 
	
	We first characterize the domain of integration occurring in the integral in $a(x_{kj})$. Let us define ${J_{kj} \defas (x_{kj}+J)} $ then
	\begin{align*}
	&\left( J^c \times J^c \cup J_{kj}^c \times J_{kj}^c \right)^c =\left( J^c \times J^c\right)^c \cap \left(J_{kj}^c \times J_{kj}^c \right)^c\\
	=&\left( (J \times J) \cup (J^c \times J) \cup (J \times J^c) \right) \cap \left((J_{kj} \times J_{kj}) \cup (J_{kj}^c \times J_{kj}) \cup (J_{kj} \times J_{kj}^c)\right) \\
	=&(C \times C) \cup (D_k \times C) \cup (C \times D_k)
	\cup (D_j \times C) \cup (J^c \cap J_{kj}^c \times C) \cup (D_j \times D_k)\\
	&\cup (C \times D_j) \cup (D_k \times J_{kj}) \cup (C \times J^c \cap J_{kj}^c),
	\end{align*}
	where we set $C \defas J \cap J_{kj}$, $D_k \defas J  \cap J_{kj}^c  $ and $D_j \defas  J^c \cap J_{kj} $. By exploiting the symmetry of the integrand we thus get
	\begin{align*}
	a_{kj}/ h ^{2d} = & \int_{J \cap J_{kj} } \int_{(J \cap J_{kj})}  F(x,y;x_{kj})dydx\\
	&+2\int_{J \cap J_{kj}} \int_{J_{kj}^c \cap J }  F(x,y;x_{kj})dydx \\
	&+2\int_{J \cap J_{kj}} \int_{J_{kj} \cap J^c }  F(x,y;x_{kj})dydx \\
	&+2\int_{J \cap J_{kj}} \int_{J^c \cap J_{kj}^c}  F(x,y;x_{kj})dydx \\
	&+2\int_{J \cap J_{kj}^c} \int_{J^c \cap J_{kj} }  F(x,y;x_{kj})dydx ,
	\end{align*}
	where $F(x,y;z) =  (\varphi(y-z)-\varphi(x-z))(\varphi(y)-\varphi(x))\g( h  y, h  x)$. Note that this representation holds for a general setting without assuming (A1) and (A2). 
	
	For implementation purpose we additionally require from now on:
	\begin{itemize}
		\item[(A3)] We assume that $R \geq diam(\Omega) = ||\bb- a ||_2$, such that $\Omega \subset B_R(x)$ for all $x \in \Omega$.
	\end{itemize}	
	This third assumption simplifies the domain of integration in the occurring integrals in the sense that we can omit the intersection with the ball $B_{R}( h  x)$. This coincides with the application to space-fractional diffusion problems where we aim to model $R\to \infty$ (see Section \ref{sec:numres}).
	Furthermore, since we can construct the whole stiffness matrix $A$ from the first row $M$, it is convenient to introduce the following $\bL^d$-dimensional vectors:
	\begin{align*}
	\sing_k &\defas { h }^{2d} \int_{J \cap J_k  } \int_{J \cap J_k}  F(x,y;E^{\bL}(k))dydx \\
	&+2{ h }^{2d}\int_{J \cap J_k } \int_{J_k^c \cap J }  F(x,y;E^{\bL}(k))dydx \\
	&+2{ h }^{2d}\int_{J \cap J_k } \int_{J_k \cap J^c }  F(x,y;E^{\bL}(k))dydx, \\
	\rad_k &\defas 2{ h }^{2d}\int_{J \cap J_k } \int_{  (J \cup J_k )^c}  F(x,y;E^{\bL}(k))\mathcal{X}_{B_R( h  x)}( h  y) dydx ,\\
	\dis_k &\defas 2{ h }^{2d}\int_{J \cap J_{k}^c} \int_{J^c \cap J_k }  F(x,y;E^{\bL}(k))dydx ,
	\end{align*}
	for $0 \leq k < \bL^d$, where $J_k \defas J_{0k} = E^{\bL}(k)+J$. Since $x_{k0} = \tfrac{x_k-x_0}{ h } = E^{\bL}(k)$, we have
	${M = \sing +\rad+ \dis}$. As we will see in the subsequent program, each of these vectors requires a different numerical handling which justifies this separation. In these premises we point out, that on the one hand we may touch possible singularities of the kernel function along the integration in $\sing_k$. On the other hand, the computation of $\rad_k$ may require the integration over a ``large'' domain if $R \to \infty$ (e.g. fractional kernel). Both are numerically demanding tasks and complicate the assembling process. In contrast to that, the computation of $\dis_k$ turns out to be numerically viable without requiring a special treatment.
	However, we fortunately find for $k$ with $J \cap J_k =  \emptyset$ that $\sing_k=0 =\rad_k$. Hence, it is worth identifying those indices and treat them differently in the assembling loop. Therefore, from $J = \dot\bigcup_{i=0}^{2d} (\square -v_i )$ we deduce that $J \cap J_k \neq \emptyset$ if and only if $k  = E^{\bL}(v_i)$ for an index $0 \leq i < 2^d$. As a consequence, we only have to compute $\sing_k$ and $\rad_k$ for $k \in \idx_0 \defas \lk E^{\bL}(v_i) \colon 0 \leq i < 2^d  \rk  $. We can cluster these indices even more. For that reason let us define on $\idx_0$ the equivalence relation
	$$k \sim j :\Leftrightarrow \exists~ \text{permutation matrix}~ P \in \R^{\bL^d \times \bL^d}\colon ~E^{-\bL}(k)  =P E^{-\bL}(j). $$
	Then due to Remark \ref{L1_radial} after Theorem \ref{multilevel_toeplitz} we have to compute the values $\sing_k$ and $\rad_k$ only for 
	$k \in  \idx_s^k\defas \lk [j]_\sim \colon j \in \idx_0 \rk$. In order to make this more precise, we figure out that the quotient set can further be specified as
	$$\idx_s^k = \lk [E^{\bL}(z)]_\sim  \colon z \in S \rk, $$
	where 
	$$S \defas \lk (0, \ldots, 0), (1,0, \ldots, 0),(1,1,0, \ldots, 0), \ldots,  (1,1, \ldots, 1)\rk \subset \lk v_i \colon 0 \leq i < 2^d\rk$$
	with $|S| = d+1 \leq 2^d$. With other words, we group those $v_i$ which are permutations of one another. The associated indices $0 \leq i < 2^d$ are thus given by
	$\lk E^{(2,\ldots,2)}(z) \colon z \in S \rk =:\idx_s^i$. 
	
\begin{figure}[h]
		\centering
		\def\svgwidth{0.5\textwidth}
		\begingroup%
		\makeatletter%
		\providecommand\color[2][]{%
			\errmessage{(Inkscape) Color is used for the text in Inkscape, but the package 'color.sty' is not loaded}%
			\renewcommand\color[2][]{}%
		}%
		\providecommand\transparent[1]{%
			\errmessage{(Inkscape) Transparency is used (non-zero) for the text in Inkscape, but the package 'transparent.sty' is not loaded}%
			\renewcommand\transparent[1]{}%
		}%
		\providecommand\rotatebox[2]{#2}%
		\ifx\svgwidth\undefined%
		\setlength{\unitlength}{125.85684468bp}%
		\ifx\svgscale\undefined%
		\relax%
		\else%
		\setlength{\unitlength}{\unitlength * \real{\svgscale}}%
		\fi%
		\else%
		\setlength{\unitlength}{\svgwidth}%
		\fi%
		\global\let\svgwidth\undefined%
		\global\let\svgscale\undefined%
		\makeatother%
		\begin{picture}(1,1.04862016)%
		\put(0,0){\includegraphics[width=\unitlength,page=1]{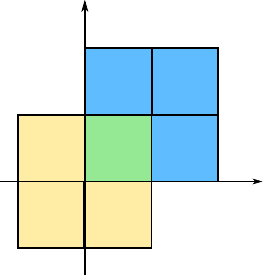}}%
		\put(0.37422495,0.46056874){\color[rgb]{0,0,0}\makebox(0,0)[lb]{\smash{$J \cap J_k$}}}%
		\put(0.12382366,0.1870326){\color[rgb]{0,0,0}\makebox(0,0)[lb]{\smash{$J \cap J_k^c$}}}%
		\put(0.63718607,0.72801492){\color[rgb]{0,0,0}\makebox(0,0)[lb]{\smash{$J^c \cap J_k$}}}%
		\put(0.27163072,0.63764145){\color[rgb]{0,0,0}\makebox(0,0)[lb]{\smash{$v_1$}}}%
		\put(0.26823475,0.38221465){\color[rgb]{0,0,0}\makebox(0,0)[lb]{\smash{$v_0$}}}%
		\put(0.52636394,0.62942166){\color[rgb]{0,0,0}\makebox(0,0)[lb]{\smash{$v_3$}}}%
		\put(0.53183302,0.3786976){\color[rgb]{0,0,0}\makebox(0,0)[lb]{\smash{$v_2$}}}%
		\put(0,0){\includegraphics[width=\unitlength,page=2]{indexsets.pdf}}%
		\put(0.0438953,0.67457862){\color[rgb]{0,0,0}\makebox(0,0)[lb]{\smash{$D_i = \lk 0 \rk$}}}%
		\put(0,0){\includegraphics[width=\unitlength,page=3]{indexsets.pdf}}%
		\put(0.62341961,0.0921992){\color[rgb]{0,0,0}\makebox(0,0)[lb]{\smash{$\lk0,1,2,3 \rk \backslash D_i = \lk 1,2,3 \rk$}}}%
		\put(0,0){\includegraphics[width=\unitlength,page=4]{indexsets.pdf}}%
		\put(0.7556215,0.94222738){\color[rgb]{0,0,0}\makebox(0,0)[lb]{\smash{$D_i^c = \lk 0,1,2 \rk$}}}%
		\end{picture}%
		\endgroup%
		\caption{\small Illustration of the index sets for $d=2$, $i=3$.}
		\label{fig:indexsets}
	\end{figure}

	In addition to the preceding considerations, we also want to partition the integration domains $J \cap J_k$, $J \cap J_k^c$ and $J^c \cap J_k$, where $k\in\idx_s^k$, into cubes $(\square - v_\nu)$, such that we can express the reference basis function $\varphi$ with the help of the element basis functions $\psi_i$. This is necessary in order to compute the integrals in a vectorized fashion and obtain an efficient implementation. Let us start with $J \cap J_k$, where $k = E^{\bL}(v_i)$ for $0 \leq i < 2^d $ such that $J_k = v_i + J$. Then we define the set
	\begin{align*}
	D_i &\defas \lk 0 \leq \mu < 2^d \colon  \square - v_\mu \in J \cap J_k \rk 
	= \lk 0 \leq \mu < 2^d \colon  \exists \kappa \colon v_\mu + v_i = v_\kappa \rk\\
	&= \lk 0 \leq \mu < 2^d \colon (v_\mu + v_i)_j < 2 ~\forall~ 0 \leq j < d \rk .
	\end{align*}
	Since $J = \dot\bigcup_{i=0}^{2^d} (\square -v_i )$, we find that $J \cap J_k = \bigcup_{\nu\in D_i}( \square - v_\nu)$. With this, we readily recognize that
	$J \cap J_k^c = \bigcup_{\nu\in\lk 0, \ldots, 2^d - 1 \rk  \backslash D_i} (\square - v_\nu)$. Similarly, one can derive a set $D_i^c$, such that
	${J^c \cap J_k = v_i + \bigcup_{\nu\in D_i^c}( \square - v_\nu)}$. Figure \ref{fig:indexsets} illustrates the latter considerations and in Table \ref{indexsets} these sets are listed for dimensions $d\in \lk 1,2,3\rk$, respectively. 
	~\\
	\begin{table}[h]
		\begin{tabular}{llllll}
			\hline\noalign{\smallskip}%
			$d$&$i \in \idx_s^i$&$k = E^L(v_i)\in \idx_s^k$ & $D_i$& $D_i^c = \lk 0\leq j < i \rk$ & $\kappa(D_i, i)$   \\
			\noalign{\smallskip}\hline\noalign{\smallskip}
			1&$0$ &0 &$(0,1)$ & $\emptyset$ &$ (0,1)$ \\
			&$1$ &1 &$0$ & $0$ &$ 1$ \\
			\noalign{\smallskip}\hline\noalign{\smallskip}
			2&$0$ &$0$&$(0,1,2,3)$ &$\emptyset$ &$ (0,1,2,3)$ \\
			&$2$ &$p_0$ &$(0,1)$ & $(0,1)$&$ (2,3)$ \\
			&$3$ & $p_0 + p_1$&$0$ & $(0,1,2)$ & $ 3$ \\
			\noalign{\smallskip}\hline\noalign{\smallskip}
			3&$0$ &$0 $ & $(0,1,2,3,4,5,6,7)$ & $\emptyset$&$ (0,1,2,3,4,5,6,7)$ \\
			&$4$ &$p_0$ & $(0,1,2,3)$ & $(0,1,2,3)$& $ (4,5,6,7)$ \\
			&$6$ &$p_0 + p_1 $ & $(0,1)$ &$(0,1,2,3,4,5)$  &$ (6,7)$ \\
			&$7$ &$p_0+p_1+p_2  $ & $0$ & $(0,1,2,3,4,5,6)$ &$ 7$ \\
			\noalign{\smallskip}\hline%
		\end{tabular}
		\caption{\small Index sets for the implementation.}\label{indexsets}
	\end{table} 
	With this at hand, we can now have a closer look at the vectors $\sing$, $\rad$ and $\dis$ and put them into a form which is suitable for the implementation.
	\subsubsection*{Vector sing}
	Let $i \in \idx_s^i$, such that $E^{\bL}(k) = v_i$. Since $J \cap J_k = \bigcup_{\nu\in D_i} \square - v_\nu$, we can transform the first integral in $\sing$ as follows
	\begin{align*}
	\int_{J \cap J_k } \int_{J \cap J_k}  F(x,y;E^{\bL}(k))dydx
	&= \sum_{\nu \in D_i} \sum_{\mu \in D_i} \int_{\square - v_\nu} \int_{\square - v_\mu}  F(x,y; v_i)dydx \\
	&= \sum_{\nu \in D_i} \sum_{\mu \in D_i} \int_{\square } \int_{\square }  F(x- v_\nu,y- v_\mu; v_i)dydx .
	\end{align*}
	By definition of $D_i$ we have for $\mu \in D_i$ that $v_\mu + v_i = v_{\kappa(\mu, i)}$ with $\kappa(\mu,i) = E^{(2,\ldots,2)}(v_i+v_\mu)$ and since $$F(x- v_\nu,y- v_\mu; v_i) =  (\varphi(y-v_{\kappa(\mu, i)})-\varphi(x-v_{\kappa(\nu, i)}))(\varphi(y-v_\mu)-\varphi(x-v_\nu))\g( h (y-v_\mu), h (x-v_\nu))$$ we find due to (\ref{reference_basisfunction}) that
	\begin{align*}
	& \sum_{\nu \in D_i} \sum_{\mu \in D_i} \int_{\square } \int_{\square }  F(x- v_\nu,y- v_\mu; v_i)dydx \\
	=& \sum_{\nu \in D_i} \sum_{\mu \in D_i} \int_{\square } \int_{\square }  ((\psi_{\kappa(\mu, i)}(y)-\psi_{\kappa(\nu, i)}(x))(\psi_{\mu}(y)-\psi_{\nu}(x)))\g( h (y- v_\mu), h (x- v_\nu))dydx.     
	\end{align*}
	We separate the case $\mu = \nu$ since the kernel may have singularities at $(x,y)$ with $x=y$ and therefore these integrals need a different numerical treatment. At this point we find that 
	\begin{align*}
	&\sum_{\nu \in D_i}\int_{\square } \int_{\square }  ((\psi_{\kappa(\nu, i)}(y)-\psi_{\kappa(\nu, i)}(x))(\psi_{\nu}(y)-\psi_{\nu}(x)))\g( h  y, h  x)dydx\\
	=&|D_i|\int_{\square } \int_{\square }  ((\psi_{i}(y)-\psi_{i}(x))(\psi_{0}(y)-\psi_{0}(x)))\g( h  y, h  x)dydx .
	\end{align*}
	This simplification follows from some straightforward transformations exploiting assumption (A2) and the fact that each element basis function $\psi_\nu$ can be expressed by $\psi_0$ through the relation
	$\psi_\nu \circ g_\nu = \psi_0 ,$ where $g_\nu(x) \defas v_\nu + R_\nu x$ for an appropriate rotation matrix $R_\nu$. With this observation we also find that the other two integrals in $\sing$ are equal, i.e.,
	\begin{align*}
	\int_{J \cap J_k} \int_{J_k^c \cap J }  F(x,y;E^{\bL}(k))dydx 
	=\int_{J \cap J_k} \int_{J_k \cap J^c }  F(x,y;E^{\bL}(k))dydx.
	\end{align*}
	By exploiting again that $J \cap J_k = \bigcup_{\nu\in D_i} \square - v_\nu$ and thus $J \cap J_k^c = \bigcup_{\nu\in\lk 0, \ldots, 2^d - 1 \rk  \backslash D_i} \square - v_\nu$ we obtain
	\begin{align*}
	&\int_{J \cap J_k} \int_{J_k^c \cap J }  F(x,y;E^{\bL}(k))dydx \\
	=&-\sum_{ \nu \in D_i} \sum_{\mu \in   \lk 0, \ldots,  2^d \rk  \backslash D_i}\int_{\square } \int_{\square } \psi_{\kappa(\nu, i)}(x) (\psi_{\mu}(y)-\psi_{\nu}(x))\g( h (v_\nu- v_\mu), h (x-y ))dydx.
	\end{align*}
	Note that by definition of $D_i$ we have that for $\mu \in   \lk 0, \ldots, 2^d \rk  \backslash D_i$ there is no $\kappa \in \lk 0, \ldots,  2^d \rk$ such that $v_\mu + v_i = v_\kappa$ and therefore $\varphi(y - (v_\mu + v_i)) = 0$ for all $y \in \square$. All in all we have
	\begin{framed}
		\begin{align*}
		&\sing_k/ h ^{2d} =\\
		& |D_i|\int_{\square } \int_{\square }  ((\psi_{i}(y)-\psi_{i}(x))(\psi_{0}(y)-\psi_{0}(x)))\g( h  y, h  x)dydx  \\
		&+ \sum_{\nu \in D_i}[ \sum_{\mu \in D_i \mu \neq \nu} \int_{\square } \int_{\square }  ((\psi_{\kappa(\mu, i)}(y)-\psi_{\kappa(\nu, i)}(x))(\psi_{\mu}(y)-\psi_{\nu}(x)))\g( h (v_\mu- v_\nu), h (y- x))dydx \\	 
		&-4 \sum_{\mu \in   \lk 0, \ldots,  2^d\rk  \backslash D_i}\int_{\square } \int_{\square } \psi_{\kappa(\nu, i)}(x) (\psi_{\mu}(y)-\psi_{\nu}(x))\g( h (v_\mu- v_\nu), h (y-x ))dydx ].
		\end{align*}
	\end{framed}
	\subsubsection*{Vector rad}
	Let $i \in \idx_s^i$, such that $E^{\bL}(k) = v_i$. Proceeding as above we obtain
	\begin{align*}
	\rad_k &= 2 h ^{2d}\int_{J \cap J_k} \int_{ (J \cup J_k)^c}  F(x,y;v_i)dydx \\
	&= 2 h ^{2d}\sum_{\nu \in D_i}\int_{\square} \psi_\nu(x)\psi_{\kappa(\nu,i)}(x)  \int_{   (J \cup J_k)^c} \g( h (x-v_\nu),  h  y) ( h  y)dydx.
	\end{align*}
	Let us define
	\begin{align*}
	P_\nu(x) \defas &  \int_{   (J \cup J_k)^c} g( h (x-v_\nu),  h  y) \mathcal{X}_{B_{R}( h (x-v_\nu))}( h  y)dy,
	\end{align*}
	where we consider $\g(x,y) = g(x,y) \mathcal{X}_{B_{R}(x)}(y)$ in accordance with remark (\ref{kernel_indicatorfunction}).
	Again, we can use $\psi_\nu \circ g_\nu = \psi_0$ in order to show by some straightforward transformations that
	$$\int_{\square} \psi_0(x)\psi_{i}(x) P_0(x)dx = \int_{\square}\psi_\nu(x)\psi_{\kappa(\nu,i)}(x)  P_\nu(x)dx $$
	for all $\nu \in D_i$. Hence, we get
	\begin{framed}
		$$\rad_k = 2  h ^{2d}|D_i|\int_{\square} \psi_0(x)\psi_{i}(x) P_0(x)dx .$$
	\end{framed}
	
	\subsubsection*{Vector dis}
	Now we distinguish between the case where $\rad_k$ and $\sing_k$ are zero and the complement case. First let $i \in \idx_s^i$, such that $E^{-\bL}(k) = v_i$, then we find
	\begin{align*}
	\dis_k &\defas 2{ h }^{2d}\int_{J_k^c \cap J} \int_{J_k \cap J^c }  F(x,y;E^{-\bL}(k))dydx \\
	& = -2{ h }^{2d} \sum_{\nu \in \lk 0,\ldots,  2^d\rk  \backslash D_i} \sum_{\mu \in  D_i^c} \int_\square \int_\square   \psi_\mu(y)\psi_\nu(x) \g( h (v_\mu - v_\nu - v_i), h (y-x)) dydx.
	\end{align*}
	Now let $k \neq E^{-\bL}(v_i)$ for any $0 \leq i < 2^d$, then $J_k^c \cap J = J$ and $J_k \cap J^c =J_k$ such that
	\begin{align*}
	\frac{\dis_k}{2{ h }^{2d}}= \int_{J} \int_{J_k}  F(x,y;E^{-\bL}(k))dydx 
	= -\int_{J} \int_{J_k}  \varphi(y+E^{-\bL}(k))\varphi(x)\g( h (y+E^{-\bL}(k)), h  x)dydx.
	\end{align*}
	Since $J_k = E^{-\bL}(k) + J $ by definition and $J = \dot\bigcup_{i=0}^{2^d} (\square -v_i )$ we obtain
	\begin{align*}
	\dis_k
	&=-2{ h }^{2d} \sum_{0 \leq \nu < 2^d} \sum_{0 \leq \mu < 2^d} \int_{\square} \int_{\square} \psi_\nu(x)\psi_\mu(y)\g( h (v_\mu-E^{-\bL}(k)-v_\nu), h (y-x))dydx.
	\end{align*}
	All in all we conclude
	\begin{framed}
		\begin{align*}
		&\dis_k/-2{ h }^{2d} \\
		= &\begin{cases}
		\sum_{\nu \in \lk 0,\ldots,  2^d\rk  \backslash D_i} \sum_{\mu \in  D_i^c} \int_\square \int_\square  \psi_\nu(x) \psi_\mu(y) \g( h (v_\mu - v_\nu - E^{-\bL}(k)), h (y-x)) dydx: &k = E^{\bL}(v_i)\\
		\sum_{0 \leq \nu < 2^d} \sum_{0 \leq \mu < 2^d} \int_{\square} \int_{\square} \psi_\nu(x)\psi_\mu(y)\g( h (v_\mu-v_\nu-E^{-\bL}(k)), h (y-x))dydx: &else.
		\end{cases}
		\end{align*}
	\end{framed}
	
	\subsubsection*{Source term}
	We compute $b^{ h } \in \R^{\bL^d} $ by
	\begin{framed}
		\begin{align*}
		b^{ h }_k = \int_{\Omega} f \hk dx
		=h^d \sum_{\nu = 0}^{2^d - 1} \int_{\square} f(x_k +  h (v- v_\nu)) \psi_{\nu}(v)dv
		=h^d 2^d \int_{\square} f(x_k  +  h  v) \psi_{0}(v)dv,
		\end{align*}
	\end{framed}
	~\\
	where the last equality follows again from considering $\psi_\nu \circ g_\nu = \psi_0$.

	\section{Solving procedure} \label{sec:solving}
	Now we discuss how to solve the discretized, fully populated multilevel Toeplitz system. The fundamental procedure uses an efficient implementation for the matrix-vector product of multilevel Toeplitz matrices, which is then delivered to the conjugate gradient (CG) method.
	
	Let us first illuminate the implementation of the matrix-vector product $Tx , $
	where $T \in \R^{\bL^d \times \bL^d}$ is a symmetric $d$-level Toeplitz matrix of order $\bL = (L_0, \ldots, L_{d-1})$ and $x$ a vector in $\R^{\bL^d}$. The crucial idea is to embed the Toeplitz matrix into a circulant matrix for which matrix-vector products can be efficiently computed with the help of the discrete fourier transform (DFT) \cite{structuredmatrices}. Here, by a $d$-level circulant matrix, we mean a matrix $C\in \R^{\bL^d \times \bL^d}$, which satisfies
	$$C_{\bi\bj} = C((\bi-\bj)\mod \bL), $$
	where $\left(\bi \mod \bL \right)\defas (i_k \mod L_k)_k$. In the real symmetric case, such that $a_{\bi\bj} = a(|\bi - \bj|)$ as it is present in our setting, $T$ can be reconstructed from its first row $R \defas (T_{0i})_i \in \R^{\bL^d}$. Since circulant matrices are special Toeplitz matrices, the same holds for these matrices as well. Due to the multilevel structure it is convenient to represent $T$ by a tensor $t$ in $\R^{L_0 \times \cdots \times L_{d-1}}$, which is composed of the values contained in $R$. More precisely we define
	$$t(\bi) \defas R_{E^{\bL}(\bi)} $$
	for $ \bi \in \prod_{i=0}^{d-1} \lk 0,\ldots, L_i-1\rk$ with $E^{\bL}$ from above. Now $t$ can be embedded into the tensor representation $c\in \R^{2L_0 \times \cdots \times 2L_{d-1}}$ of the associated $d$-level circulant matrix by
	$$c(\bi) \defas t(\hat{i}_0, \ldots, \hat{i}_{d-1}), $$
	where
	$$\hat{i}_k \defas \begin{cases}
	i_k: &i_k<L_k, \\
	0: &i_k = L_k, \\
	2L_k-i_k: &else.
	\end{cases}$$
	We note that 
	$t = c([0:L_0-1], \ldots, [0:L_{d-1}-1]) $. Thus, we can use Algorithm \ref{product} to compute the product $T \cdot x$, where the DFT is carried out by the fast fourier transform (FFT).
	\begin{algorithm}
		\caption{Matrix-vector product for multilevel Toeplitz matrices}
		\label{product}
		~\\
		INPUT: 		$t \in \R^{L_0 \times \cdots \times L_{d-1}}$ representing $T \in \R^{\bL^d \times \bL^d}$, $x \in \R^{\bL^d} $\\
		OUTPUT: 	$y = Tx$\\
		\\
		1. Construct $c\in \R^{2L_0 \times \cdots \times 2L_{d-1}}$ by $c(\bi) \defas t((\hat{i}_0, \ldots, \hat{i}_{d-1}) $\\
		2. Construct $x' \in \R^{2L_0 \times \cdots \times 2L_{d-1}}$ by $$x'(\bi)\defas \begin{cases}
		x_{E^L(\bi)}: & \bi \in \prod_{i=0}^{d-1} \lk 0,\ldots, L_i-1\rk \\
		0: else
		\end{cases}$$\\
		3. Compute $\Lambda = FFT_L(c) $\\
		4. Compute $z = FFT_L(x')$\\
		5. Compute $w = \Lambda  z $ (pointwise)\\
		6. Compute $y' = FFT^{-1}_L(w)$\\
		7. Construct $y \in \R^{L_0 \times \cdots \times L_{d-1}}$ by $y(\bi) = y'(\bi)$ for $ \bi \in \prod_{i=0}^{d-1} \lk 0,\ldots, L_i-1\rk $\\
		8. Return $y.reshape(\bL^d)$
	\end{algorithm}
	In the Python code we use the library pyFFTW (\url{https://hgomersall.github.io/pyFFTW/}) to perform a parallelized multidimensional DFT, which is a pythonic wrapper around the C subroutine library FFTW (\url{http://www.fftw.org/}). Furthermore, we want to note, that an MPI implementation for solving multilevel Toeplitz systems in this fashion is presented in \cite{MPIproduct}, which inspired us to apply the upper procedure.\\ Finally, with this algorithm at hand, we employ a CG method, as it can be found for example in \cite{nocedale:numopt}, to obtain the solution of the discretized system (\ref{discretizedproblem}).

	\section{Numerical Experiments}\label{sec:numres}
	In this last section we want to complete the previous considerations by presenting numerical results in 1d, 2d and for the first time also in 3d. We now specify the nonlocal diffusion operator $-\Lu$ by choosing the truncated fractional kernel $\g(x,y) = \frac{c_{d,s}}{2||y-x||^{d+2s}}\mathcal{X}_{B_R(x)}(y)$ and shortly recall how the fractional Laplace operator $(-\Delta)^s$ crystallizes out as special case of $-\Lu$. Furthermore we describe in detail the numerical integration and finally discuss the results of the implementation.
	
	\subsection{Relation to space-fractional diffusion problems}
	We follow \cite{fractlapl} and define the action of the \emph{fractional Laplace operator $(-\Delta)^s$} on a function ${u:\R^d \to \R}$ by
	\begin{align*}
	(-\Delta)^s u(x) \defas c_{d,s} \int_{\R^d} \frac{u(x)- u(y)}{||y-x||^{d+2s}}dy,
	\end{align*}
	where
	$c_{d,s} \defas s 2^{2s} \frac{\Gamma(\tfrac{d+2}{2})}{\Gamma(\tfrac{1}{2}) \Gamma(1-s)}  .$
	The homogeneous steady-state \emph{space-fractional diffusion problem} then reads as 
	\begin{align}
	(-\Delta)^s u(x) &= f(x), ~~~x \in \Omega, \nonumber\\
	u(x) &= 0, ~~~~~~~ x \in \Omega^c \label{fractprob}.
	\end{align}
	Thus, choosing the truncated fractional kernel $\g(x,y) \defas \frac{c_{d,s}}{2||y-x||^{d+2s}}\mathcal{X}_{B_R(x)}(y)$, which satisfies the conditions (\ref{case_1}) and (\ref{gen_assumption}), the nonlocal Dirichlet problem (\ref{prob}) can be considered as a truncated version of problem (\ref{fractprob}). In \cite{fractlapl} the authors show that the weak solution $u_R$ of the truncated problem (\ref{prob}), for some interaction radius $R>0$, converges to the weak solution $u_\infty$ of (\ref{fractprob}) as $R \to \infty$. We recall the corresponding result for completeness.
	\begin{prop}[{\cite[Theorem 3.1]{fractlapl}}]
		Let $u_R \in V_c(\OuO)$ and $u_\infty \in H_\Omega^s(\R^d) \defas \lk u \in H^s(\R^d) \colon u_{|\Omega^c} \equiv 0 \rk$ denote the weak solutions of (\ref{prob}) and (\ref{fractprob}) respectively. Then $$||u_\infty - u_R ||_{H^s(\OuO)} \leq \frac{K_d}{C_1^2 s(R-I)^{2s}}||u_\infty||_{L^2(\Omega)},$$ where $I \defas \min\lk L\in \R \colon \Omega \subset B_L(x)~~ \forall ~x \in \Omega \rk$, $C_1$ is the equivalence constant from (\ref{space_eq}) and $K_d$ is a constant depending only on the space dimension $d$.
	\end{prop}
	
	\subsection{Numerical computation of the integrals} 
	In this subsection we want to point out how we numerically handle the occurring integrals. 
	\subsubsection{Nonsingular integrals}
	We mainly have to compute integrals of the form
	$
	\int_{\square} \int_{\square} g(x,y) dydx \label{referenceintegral},
	$
	where $\square = [0,1]^d$ and $g \colon \R^d \times \R^d \to \R$ is a (typically smooth) function, which we assume to have no singularities in the domain $\square \times \square$. We approximate the value of this integral by employing a $n$-point Gauss-Legendre quadrature rule in each dimension. More precisely,  we built a $d$-dimensional tensor grid $\X \in \R^{d \times n^{d}}$ with associated weights $\texttt{W}^{single} \in \R^{n^d}$, such that $\int_{\square} g(x,y) dy \approx \sum_{i=0}^{n^{d}-1} g(x, \X_i)  \texttt{W}^{single}_i$ for $x \in \X$.
	Finally, we define the arrays
	\begin{align*}
	\texttt{V} &\defas (\X, \X, \ldots, \X) \in \R^{d \times n^{2d}}, \\
	\texttt{Q} &\defas (\X_0, \ldots, \X_0,\X_1, \ldots, \X_1, \ldots , \X_{n^d - 1}, \ldots, \X_{n^d - 1} ) \in \R^{d \times n^{2d}} \
	\end{align*}
	with associated weights $\texttt{W}^{double}\in \R^{n^{2d}}$, such that we finally arrive at the following quadrature rule:
	\begin{align*}
	\int_{\square} \int_{\square} g(x,y) dydx \approx \sum_{i=0}^{n^{2d}-1} g(\texttt{V}_i, \texttt{Q}_i)  \texttt{W}^{double}_i = \left(g(\texttt{V},\texttt{Q}) \cdot \texttt{W}^{double}\right).sum().
	\end{align*}
	\subsubsection{Singular integrals}
	As we have seen for the space-fractional diffusion problem, the kernel function may come along with singularities at $(x,y)$ with $x=y$. Therefore we start with a general observation, which paves the way for numerically handling these singularities. 
	Let $f \colon \R^d \times \R^d \to \R$ be a symmetric function, i.e., $f(x,y) = f(y,x)$, and let us further define the sets
	$M \defas \lk (x,y) \in \square\times \square \colon y_d \in [0,x_d] \rk$
	and
	$M' \defas \lk (x,y) \in \square\times \square \colon (y,x) \in M\rk.$
	Then it is straightforward to show that $M \cup M' = \square \times \square$ and
	\begin{align*}
	M \cap M' = &\lk (x,y)\in \square \times \square \colon x_d = y_d \rk \\
	= &\lk (v^{d-1},z,w^{d-1},z) \in \R^{2d} \colon (v^{d-1},w^{d-1},z) \in [0,1]^{2d-1} \rk,
	\end{align*}
	such that $\lambda_{2d}(M \cap M')= 0$. Hence, we find that
	\begin{align*}
	\int_{\square \times \square} f d\lambda_{2d} 	
	&= \int_{M} f d\lambda_{2d} + \int_{M'} f d\lambda_{2d} - \int_{M \cap M'} f d\lambda_{2d}\\
	&=\int_{M} f d\lambda_{2d} + \int_{M'} f d\lambda_{2d} .
	\end{align*}
	From the symmetry of $f$ we additionally deduce that
	$\int_{M} f d\lambda_{2d} = \int_{M'} f d\lambda_{2d} $
	and therefore we finally obtain
	\begin{align*}
	\int_{\square \times \square} f d\lambda_{2d} 	= 2 \int_{M} f d\lambda_{2d}.
	\end{align*}
	This observation can now be applied to the singular integrals occurring in the vector $\sing$ such that
	\begin{align*}
	&\int_{\square}\int_{\square}((\psi_{i}(y)-\psi_{i}(x))(\psi_{0}(y)-\psi_{0}(x)))\g( h  y, h  x)dydx   \\
	= &2\int_{[0,1]^d}\int_{[0,1]^{d-1}\times [0,x_{d-1}]} ((\psi_{i}(y)-\psi_{i}(x))(\psi_{0}(y)-\psi_{0}(x)))\g( h  y, h  x)dydx .
	\end{align*}
	The essential advantage of this representation relies on the fact that the singularities are now located on the boundary of the integration domain. Thus, we do not evaluate the integrand on its singularities while using quadrature points which lie in the interior. We extend the one-dimensional adaptive (G7,K15)-Gauss-Kronrod quadrature rule to $d$-dimensional integrals by again tensorising the one-dimensional quadrature points. Moreover, in order to take full advantage of the Gauss-Kronrod quadrature, we divide the set $[0,1]^{d-1}\times [0,x_{d-1}]$ into $2^{d-1}$ disjoint rectangular subsets such that the singularity $x$ is located at a vertex. The latter partitioning reinforces the adaptivity property of the Gauss-Kronrod quadrature rule.
	
	\subsubsection{Integrals with large interaction horizon}
	Now we discuss the quadrature of
	\begin{align*}
	P_0(x) =~ &\int_{ (J \cup J_k)^c} g( h  x,  h  y)\mathcal{X}_{B_{R}( h  x)}( h  y) dy \\
	=~&(1/h^d)\int_{ (I_0 \cup I_k )^c} g( h  x, y-( a + h  ))\mathcal{X}_{B_{R}( h  x)}(y-( a + h  )) dy.
	\end{align*}
	Recall that we consider $\g(x,y) = g(x,y) \mathcal{X}_{B_{R}(x)}(y)$ as in (\ref{kernel_indicatorfunction}). We carry out two simplifications for the implementation, which are mainly motivated by the fact that $\g(x,y) \to 0$ as $y  \to \infty$ for kernels such as the fractional one. First, since $h \to 0$, we set $B_{R}( h  x) \equiv B_{R}(0)$. Especially when $R$ is large and $h$ small, this simplification does not significantly affect the value of the integral. Second, we employ the $||\cdot||_\infty$-norm for the ball $B_R(0)$ instead of the $||\cdot||_2$-norm. Hereby we also loose accuracy in the numerical integration but it simplifies the domain of integration in the sense that we can use our quadrature rules for rectangular elements. The latter simplification can additionally be justified by the fact that we want to model $R \to \infty$ for the fractional kernel and therefore have to truncate the $||\cdot||_2$-ball in any case. Consequently, we are concerned with the quadrature of the integral
	\begin{align*}
	\int_{  B_{R}^{||\cdot||_\infty}( a  +  h )   \backslash (I_0 \cup I_k )} g( h  x, y-(a+ h )) dy
	\end{align*}
	for $x \in \X$. For this purpose, we define the box $\mathcal{B}\defas\prod_{i=0}^{d-1} [a_i-\lambda,a_i+\lambda]$ for a constant $R\geq \lambda > 2h$ such that $(I_0 \cup I_k ) \subset \mathcal{B}$ and we partition
	$$B_{R}^{||\cdot||_\infty}( a  +  h )  = \mathcal{B} \cup B_{R}( a  +  h ) \backslash \mathcal{B} .$$
	Thus, we obtain
	\begin{align*}
	P_0(x)h^d&\approx\int_{  \mathcal{B}  \backslash (I_0 \cup I_k )} g( h  x, y-( a  +  h )) dy + \int_{  B_{R}^{||\cdot||_\infty}( a  +  h )  \backslash  \mathcal{B}} g( h  x, y-( a  +  h )) dy.
	\end{align*}
	We discretize $ \mathcal{B}$ with the same elements which we used for $\Omega$. This is convenient for two reasons. On the one hand we capture the critical values of the kernel, which in case of singular kernels typically decrease as $y\to \infty$. On the other hand, we have to leave out the integration over $(I_0 \cup I_k )$, which then can easily be implemented since we use the same discretization. However, this results in $\bN_2^d$, where $N_2^i \defas 2\lambda /  h $, hypercubes with base points 
	$ y_j \defas ( a -\lambda e) +   h  E^{-\bN_2}(j)$ such that
	\begin{align*}
	\int_{   \mathcal{B}  \backslash (I_0 \cup I_k )} g( h  x, y) dy&= \sum_{j=0,j \notin R_i}^{\bN_2^d -1} \int_{ y_j +  h \square} g( h  x,  y-( a  +  h )) dy\\
	&= h^d \sum_{j=0,j \notin R_i}^{\bN_2^d - 1} \int_{ \square} g(\lambda e + h (e-E^{-\bN_2}(j)),  h (y-x)) dy,
	\end{align*}
	where $R_i$ contains the indices for those elements, which are contained in $(I_0 \cup I_k )$. This set can be characterized as $R_i \defas \lk 0\leq j < \bN_2^d - 1 \colon  y_j +  h \square \subset  (I_0 \cup I_k )\rk$. Note that by definition we have $I_0 =  a + h (e + \bigcup_{0 \leq \nu < 2^d} \square - v_\nu)$ and since $E^{-\bL}(k) = v_i$ for $k \in \idx_s^k$ such that $x_k =  a + h (e +v_i)$ we know that $I_0^c \cap I_k =  a + h (e +v_i + \bigcup_{\nu \in   D_i^c} \square - v_\nu)$. Hence, $ y_j +  h \square \subset  (I_0 \cup I_k )$ if and only if $ y_j =  a + h (e - v_\nu)$ for $0 \leq \nu < 2^d$ or $ y_j =  a +   h (e +v_i- v_\nu)$ for $\nu \in   D_i^c$. Since $ y_j =  a -\lambda e + h  E^{-\bN_2}(j) $ we find by equating $ y_j$ with these requirements that
	\begin{align*}
	R_i =\lk E^{\bN_2}(e -  v_\nu +\lambda  h ^{-1} e) : 0 \leq \nu < 2^d \rk 
	\cup \lk E^{\bN_2}(e +v_i -  v_\nu +\lambda  h ^{-1} e): \nu \in   D_i^c \rk.
	\end{align*}
	Now we discuss the quadrature of the second integral
	$\int_{  B_{R}^{||\cdot||_\infty}( a  +  h )  \backslash  \mathcal{B}} g( h  x, y-( a  +  h )) dy$
	for $x \in \X$. Since we want to apply the algorithm to fractional diffusion, we have to get around the computational costs that occur when $R$ is large. In order to alleviate those costs we follow the idea in \cite{fractlapl} and apply a coarsening rule to discretize the domain $ B_{R}^{||\cdot||_\infty}( a  +  h ) \backslash \mathcal{B}$. Assume we have a procedure which outputs a triangulation $(z, \hat{ h })$ of $  B_{R}^{||\cdot||_\infty}( a  +  h )  \backslash \mathcal{B}$ consisting of $ N_3 $ hyperrectangles with base points $ z_j \in \R^d$ and sides of length $ \hat{ h }_j   \in \R^d$. Then we obtain
	\begin{align*}
	\int_{  B_{R}^{||\cdot||_\infty}( a  +  h )  \backslash \mathcal{B}} g( h  x, y-( a  +  h )) dy &= \sum_{j=0}^{ N_3-1 } \int_{ z_j +  \hat{ h }_j \square} g( h  x,  y-( a  +  h ) ) dy\\
	&= \sum_{j=0}^{ N_3 -1} \hat{ h }_j ^d\int_{ \square} g(- z_j+ a  +  h ,  \hat{ h }_j y - h  x ) dy.
	\end{align*} 
	All in all we thus have
	\begin{framed}
		\begin{align*}
		\rad_k  &\approx 2|D_i|  { h }^{2d}\sum_{j=0,j \notin R_i}^{\bN_2^d-1} \int_{\square}\int_{ \square}\psi_0(x)\psi_{i}(x)g(\lambda e + h (e -E^{-\bN_2}(j)),  h (y-x)) dydx \\
		&+2|D_i| h ^{d} \sum_{j=0}^{ N_3-1 } \hat{ h }_j ^d\int_{\square}\int_{ \square}\psi_0(x)\psi_{i}(x)g(- z_j+ a  +  h ,  \hat{ h }_j  y- h  x) dydx.
		\end{align*}
	\end{framed}
	This leaves space for discussion concerning the choice of an optimal coarsening rule. We use the following simple approach in our code: We decompose $  B_{R}^{||\cdot||_\infty}( a  +  h )  \backslash \mathcal{B}$ into $(3^d-1)$ $d$-dimensional hyperrectangles surrounding the box $\mathcal{B}$. Then we build a tensor grid by employing in each dimension the coarsening strategy
	$ v + i^qh_{min} $ where $v$ is a vertex of $\mathcal{B}$, $q \geq 1$ the coarsening parameter and $h_{min}$ a minimum grid size. By concatenating all arrays we obtain a triangulation $(z, \hat{ h })$ of $  B_{R}^{||\cdot||_\infty}( a  +  h )  \backslash \mathcal{B}$.

	\subsection{Numerical results}
	The implementation has been carried out in Python and the examples were run on a HP Workstation Z240 MT J9C17ET with Intel Core i7-6700 - 4 x 3.40GHz. Since we started from an arbitrary dimension throughout the whole analyzes, the codes for each dimension $d \in \lk 1,2,3\rk$ own the same structure. Depending on the dimension, one only has to adapt the index sets $\idx_s^i$, $\idx_s^k$, $D_i$ and $D_i^c$, the implementation of the quadrature rules discussed above and the $2^d$ element basis functions. The rest can be implemented in a generic way. In addition to that, we framed above the relevant representations for implementing the assembly of the first row. The implementation of the solving procedure only consists of delivering Algorithm \ref{product} to the CG method. Moreover, the codes are parallelized over 8 threads on the four Intel cores. Within the assembling process of the first row, we first compute the more challenging $(d+1)$ entries $M_k = \sing_k + \rad_k + \dis_k$ for $k \in \idx_s^k$. Here we parallelize the computations of the integrals in $\rad_k$ over the base points $y_j$ for the box and $z_j$ for the coarsening strategy. For the remaining $\bL^d - (d+1)$ indices we have that $M_k = \dis_k$ and we can simply parallelize the loop over $\lk 0\leq k < \bL^d \rk \backslash \idx_s^k$. As mentioned above, the solving process is parallelized via the parallel fourier transform pyFFTW.
	
	In all examples we consider $\Omega = [0,1]^d$ and use the truncated fractional kernel $$\g(x,y) = \frac{c_{d,s}}{2||y-x||_2^{d+2s}}\mathcal{X}_{B_R(x)}(y)$$ for $s=0.4$ and $R=T+\lambda$ where $T=2^{10}$ with coarsening parameter $q=1.5$, minimum grid size $h_{min}=10^{-2}$ and a parameter $\lambda > 2h$ for the box $\mathcal{B}$. Furthermore we consider a constant source term $f\equiv 1$ in \eqref{prob}.
	The CG method stops if a sufficient decrease of the residual $||Ax_k-b||/||b|| < 10^{-12}$ is reached.  
	We present numerical examples for $d \in \lk 1,2,3\rk$. For each grid size $ h $  we report on the number of grid points (``dofs'') and the number of CG iterations (``cg its'') as well as the CPU time (``CPU solving'') needed for solving the discretized system. 
	Furthermore we compute the energy error $||u_R^h - u_\infty||_{H^s(\OuO)}$, where $u_\infty$ is 
	a numerical surrogate taken to be the finite element solution on the finest grid, and the rate of convergence. 
	
	\subsubsection*{1d Example}
\begin{figure}[h]
		\label{fig:1d}
		\centering
		\includegraphics[width =0.49 \textwidth]{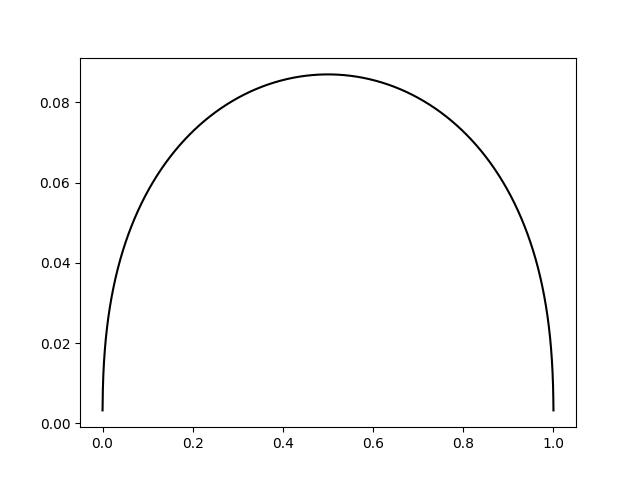}
		\caption{\small Plot of the 1d finite element solution $u^{ h }$.}	
	\end{figure}

	For the 1d example we choose $\lambda = 5$ as parameter for the box $\mathcal{B}$ and $n=7$ Gauss points for the unit interval $[0,1]$. The results are presented in Figure \ref{fig:1d} and Table \ref{tab:1d}.
	\begin{table}[h]
		
		\centering
		\begin{tabular}{llllll}
			\hline\noalign{\smallskip}
			$h$ 		& dofs 		& cg its    & energy error &rate 		&CPU solving [s]\\ 
			\noalign{\smallskip}\hline\noalign{\smallskip}
			
			$2^{-6}$ 	& 63		&16 		&2.43e-02 		&   0.50	&-	\\ 
			$2^{-7}$ 	& 127 	&24 		&1.72e-02  		& 0.51   	&0.005  		\\ 
			$2^{-8}$ 	& 255  	&34		&1.21e-02  	&  0.51		& 0.011 	\\ 
			$2^{-9}$  & 511   	& 46 	&8.47e-03  		&  0.52		& 0.029 \\ 
			$2^{-14}$ & 16,383  	& 191 	&  - 			&  - 		& 5.26 \\ 
			\noalign{\smallskip}\hline
		\end{tabular}	
		\caption{\small Results of the 1d test case. 
		}\label{tab:1d}
	\end{table}

	\subsubsection*{2d Example}
\begin{figure}[H]
		\label{fig:2d}
		\centering
		\includegraphics[width =0.49 \textwidth]{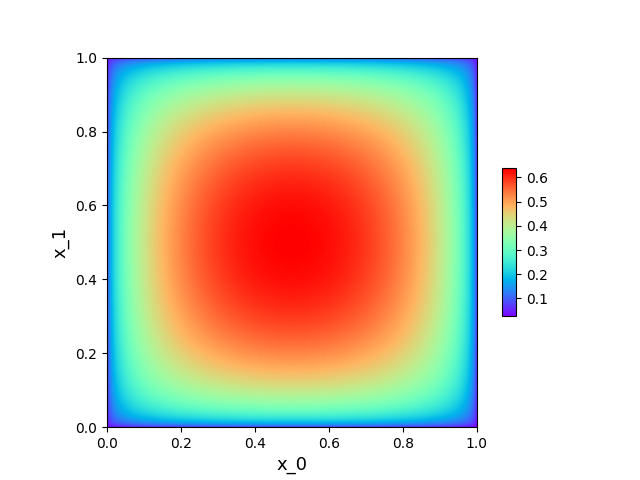}
		\caption{\small Contour plot of the 2d finite element solution $u^{ h }$.}
		
	\end{figure}

	For the 2d example we choose $\lambda = 1$ as parameter for the box $\mathcal{B}$ and $n=6$, i.e., $36$ quadrature points for the unit square $[0,1]^2$. The results are presented in Figure \ref{fig:2d} and Table \ref{tab:2d}.

	\begin{table}[H]
		
		\centering
		\begin{tabular}{ llllll}
			\hline\noalign{\smallskip}
			$h$ 		& dofs 			& cg its& energy error  &rate 	&CPU solving [s]	\\
			\noalign{\smallskip}\hline\noalign{\smallskip}
			$2^{-2}$ 	& 9 		&3   	&  	3.11e-01			& 0.50		&-   				\\
			$2^{-3}$ 	& 49  		&10   	&	2.17e-01 			& 0.51		&-     				\\
			$2^{-4}$ 	& 225 		&16   	&   1.53e-01			& 0.52		&0.01    				\\
			$2^{-5}$   & 961  		&20   	&  	1.06e-01			& 0.53		&0.03					\\
			$2^{-9}$   & 261,121  	& 58 	& 	-					& -		& 	3.15				\\
			\noalign{\smallskip}\hline
		\end{tabular}	
		\caption{\small Results of the 2d test case.}
		\label{tab:2d}
	\end{table}
	
	\subsubsection*{3d Example}
\begin{figure}[H]
		\label{fig:3d}
		\centering
		\includegraphics[width =0.5 \textwidth]{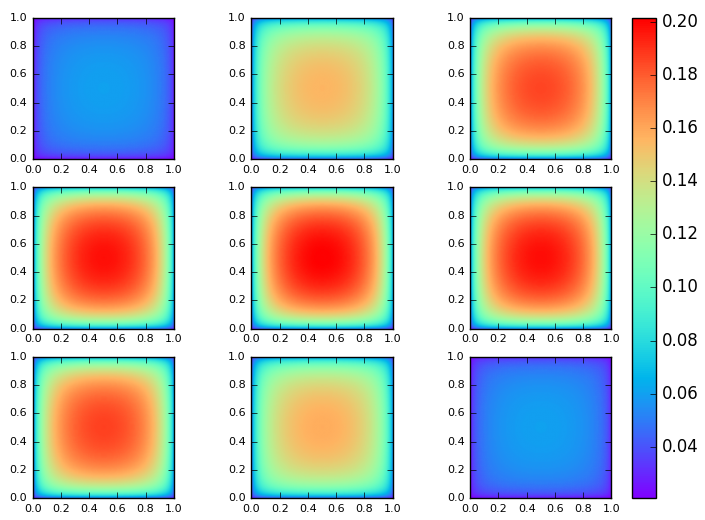}
		\caption{\small Plot of the finite element solution $u^{ h }$. In order to illustrate the 3d solution we cut the domain  $\Omega = [0,1]^3$ into nine slices along the third dimension ordered by increasing $x_2$-dimension.}
		
	\end{figure}

	For the 3d example we choose $\lambda = 0.5$ as parameter for the box $\mathcal{B}$ and $n=4$, i.e., $64$ quadrature points for the unit cube $[0,1]^3$. The results are presented in Figure \ref{fig:3d} and Table \ref{tab:3d}.
	
	\begin{table}[H]
		\centering
		\begin{tabular}{llllll}
			\hline\noalign{\smallskip}
			$h$ 		& dofs & cg its           & energy error &rate &CPU solving [min]\\
			\noalign{\smallskip}\hline\noalign{\smallskip}
			$2^{-3}$ 	& 343  &19 						& 1.37e-01 		&0.51	 			&-		\\ 
			$2^{-4}$ 	& 3,375  &20					& 9.49e-02 		&0.52	 			&-		\\ 
			$2^{-5}$ 	& 29,791  &21 					& 6.54e-02 		&0.54	 			&0.01		\\ 
			$2^{-6}$  & 250,047  & 23 					& 4.44e-02 		& 0.59				&0.03 		\\
			$2^{-9}$ & 133,432,831  & 55 				& -		&- 				&36.28		\\
			\noalign{\smallskip}\hline
		\end{tabular}	
		\caption{\small Results of the 3d test case.}
		\label{tab:3d}
	\end{table}
	
	\subsection{Discussion}
	Let us first comment on the convergence rate. Since we set $f\equiv 1$, we find that our 1d results confirm the theoretical result given in (\ref{convergence_rate}) where $\alpha = \min \lk s+r,1/2-\varepsilon \rk  \approx 0.5$. Due to the numerical results for $d>1$ one may conjecture, that this convergence result also holds for domains with less smooth boundary, such as hyperrectangles. The latter has already been shown for finite element approximations of the untruncated problem (\ref{fractprob}); see \cite[Thoerem 4.7]{borthagaray_regularity}. 

	Concerning the number of CG iterations we point out two observations. First, we generally observe for various parameters that the number of CG iterations increases as we emphasize the singularity, i.e., as $s \to 1$. This  coincides with Theorem 6.3 in \cite{wellposedness}, stating that for shape-regular and quasi-uniform meshes the following estimate for the condition number $cond(A)$ holds: $$cond(A)\leq ch^{-2s},$$ where $c>0$ is a generic constant. 
	A second observation in this regard is the surprisingly low number of CG iterations needed for the 3d case. Running the code for different dimensions, but for a fixed comparable parameter setting, where we set the domain to be the unit hypercube $[0,1]^d$ and the source term to be $f \equiv 1$, we find that for problems of the same size, i.e., with the same number of degrees of freedom, the number of CG iterations decreases as the dimension increases. Specifically, let $s=0.4$, $\lambda = 0.5$ and $T=2^{10}$. Then in Table \ref{tab:compare_dofs} we find the results where the size of the discretized system is fixed to $dofs = \bL^d \approx 250000$.
	\begin{table}[H]	
		\centering
		\begin{tabular}{ lll}
			\hline\noalign{\smallskip}
			$d$ 				& h				 	& cg its           			\\
			\noalign{\smallskip}\hline\noalign{\smallskip}
			$1$ 				& $1./250,048 \approx 2^{-18} $ 		&615				\\ 
			$2$  				& $2^{-9}$		& 58 						\\
			$3$  			  	& $2^{-6} $				&23		\\
			\noalign{\smallskip}\hline
		\end{tabular}	
		\caption{\small Results for a fixed size $\bL^d \approx 250000$, but different grid sizes $h$.}
		\label{tab:compare_dofs}
	\end{table}
	This might be due to the additional Toeplitz levels affecting the condition number of the stiffness matrix. In contrast to that, while fixing the grid size, the number of CG iterations seems to vary less comparing the different dimensions. In Table \ref{tab:compare} the reader finds the results for a fixed grid size $h= 2^{-6}$.
	\begin{table}[H]	
		\centering
		\begin{tabular}{ lll}
			\hline\noalign{\smallskip} 
			$d$ 				& dofs			 	& cg its           			\\
			\noalign{\smallskip}\hline\noalign{\smallskip}
			$1$ 				&  63		&16				\\ 
			$2$  				& 3,969	& 23 						\\
			$3$  			  	& 250,047 			&23		\\
			\noalign{\smallskip}\hline
		\end{tabular}	
		\caption{\small Results for a fixed grid size $h = 2^{-6}$, but different system sizes $dofs = \bL^d$.}
		\label{tab:compare}
	\end{table}
	However, this has to be analyzed more concretely and a thorough investigation is not the intention at this point.
	
	
	Finally, we also note, that the library pyFFTW needs a lot of memory for building the FFT object, such that we had to move the 3d computations for the finest grid to a machine with a larger RAM. One can circumvent this problem by using the sequential FFT implementation available in the NumPy library. 
	
	\section{Concluding remarks}
	We presented a finite element implementation for the steady-state nonlocal Dirichlet problem with homogeneous volume constraints on an arbitrary $d$-dimensional hyperrectangle and for translation and reflection invariant kernel functions. We use a continuous Galerkin method with multilinear element basis functions and theoretically back up our numerics with the framework for nonlocal diffusion developed by Gunzburger et al. \cite{wellposedness}. The key result showing the multilevel Toeplitz structure of the stiffness matrix is proven for arbitrary dimension and paves the way for the first 3d implementations in this area. Furthermore, we comprehensively analyze the entries of the stiffness matrix and derive representations which can be efficiently implemented. Since throughout the whole analysis we start from an arbitrary dimension, one can almost generically implement the code by adapting the implementation of the quadrature rules, the element basis functions as well as the index sets.
	
	An important extension of this work is to incorporate the case where the interaction horizon is smaller than the diameter of the domain. This complicates the integration and with that the assembling procedure, but the stiffness matrix is no more fully populated and its structure still remains multilevel Toeplitz. Having that, one can model the transition to local diffusion and access a greater range of kernels. The resulting code for 2d would then present a fast and efficient implementation, which could be used for example in image processing.
	Also, since certain kernel functions allow for solutions with jump discontinuities, also discontinuous Galerkin methods are conforming \cite{wellposedness}. In this case, one has to carefully analyze the structure of the resulting stiffness matrix, which might differ from a multilevel Toeplitz one. Moreover, an aspect concerning the solving procedure, which is not examined above, is that of an efficient preconditioner for the discretized Galerkin system. A multigrid method might be a reasonable candidate due to the simple structure of the grid (see also \cite{chen2016convergence,PANG2012693}). In general, a lot of effort has been put in the research of preconditioning structured matrices (see e.g., \cite{olshans,serra2007toeplitz,chan1988optimal,Capizzano2000AnyCP}). Since we observe a moderate number of CG iterations in our numerical examples, a preconditioner has not been implemented yet.
	
	The main drawback of our approach relies on the fact that the code is strictly limited to regular grids and is thus not applicable to more complicated domains. It is crucial that each element has the same geometry in order to achieve the multilevel Toeplitz structure of the stiffness matrix; meaning that only rectangular domains are reasonable.  
	However, one could think of a coupling strategy, which allows us to decompose a general domain into rectangular parts and the remaining parts. This is beyond the scope of this paper and left to future work.
	In contrast to that, the restriction to translation and reflection invariant kernels appears to be rather weak, since a lot of kernels treated in literature are even radial. 
	

	%
	%

	\section*{Acknowledgements}
		The first author has been supported by the German Research Foundation
		(DFG) within the Research Training Group 2126: ``Algorithmic Optimization''.

	
	\pagestyle{plain}
	
	\bibliographystyle{plain}
	\bibliography{literatur.bib}   

\end{document}